\documentclass[11pt]{amsart} \textwidth=14.5cm \oddsidemargin=1cm \evensidemargin=1cm \usepackage{amsmath}
\usepackage{amsxtra} \usepackage{amscd} \usepackage{amsthm} \usepackage{amsfonts} \usepackage{amssymb}
\usepackage{amscd}
\usepackage{eucal}
\usepackage[matrix,arrow,curve]{xy}
\usepackage[dvipsnames]{xcolor}
\newcommand{\Red}[1]{{\color{red}#1}}

\numberwithin{equation}{section} \newtheorem{thm}{Theorem} \newtheorem{cor}[thm]{Corollary}
\newtheorem{lem}[thm]{Lemma} \newtheorem{prop}[thm]{Proposition} 

\theoremstyle{definition} 

\theoremstyle{remark} \newtheorem{rem}[thm]{Remark} \newtheorem*{rem*}{Remark}

\numberwithin{thm}{section}

\theoremstyle{definition}

\theoremstyle{remark}

\newcommand{\nc}{\newcommand} \nc{\renc}{\renewcommand} \nc{\ssec}{\subsection} \nc{\sssec}{\subsubsection}
\nc{\on}{\operatorname} \nc\wt{\widetilde}

\newcommand{\iso}{\stackrel{\sim}{\longrightarrow}} 
\newcommand{\CalD}{{\mathcal D}} 

\newcommand{\Lotimes}{{\stackrel{L}{\otimes}}}

\nc{\BA}{{\mathbb{A}}}
 \nc{\BC}{{\mathbb{C}}}
\nc{\BD}{{\mathbb{D}}}
\nc{\BF}{{\mathbb{F}}} \nc{\BG}{{\mathbb{G}}} \nc{\Gmo}{{\mathbb{G}_m^{(1)}}}
\nc{\BH}{{\mathbb{H}}}
\nc{\BM}{{\mathbb{M}}} \nc{\BN}{{\mathbb{N}}} \nc{\BP}{{\mathbb{P}}}
\nc{\BQ}{{\mathbb{Q}}}
\nc{\BR}{{\mathbb{R}}} \nc{\BZ}{{\mathbb{Z}}} \nc{\BS}{{\mathbb{S}}}

\nc{\Zet}{{\mathbb{Z}}} \nc{\Ce}{{\mathbb{C}}} \nc{\Fq}{{\mathbb{F}_q}}

\nc{\CA}{{\mathcal{A}}} \nc{\A}{{\mathcal{A}}} \nc{\CB}{{\mathcal{B}}} \nc{\B}{{\mathcal{B}}}
\nc{\CC}{{\mathcal{C}}} 
\nc{\CE}{{\mathcal{E}}} \nc{\CF}{{\mathcal{F}}}
\nc{\F}{{\mathcal{F}}} \nc{\CG}{{\mathcal{G}}} \nc{\G}{{\mathcal{G}}}
\nc{\CH}{{\mathcal H}} \nc{\CK}{{\mathcal{K}}}
\nc{\CI}{{\mathcal{I}}} \nc{\CJ}{{\mathcal J}} \nc{\CL}{{\mathcal{L}}} \nc{\CM}{{\mathcal{M}}}
\nc{\CMM}{{\mathcal{M}^{\operatorname{gen}}_\hbar(-\rho)}} \nc{\CN}{{\mathcal{N}}}
\nc{\CO}{{\mathcal{O}}}
\nc{\CP}{{\mathcal{P}}} \nc{\CQ}{{\mathcal{Q}}} \nc{\CR}{{\mathcal{R}}} \nc{\CS}{{\mathcal{S}}}
\nc{\CT}{{\mathcal{T}}} \nc{\CU}{{\mathcal{U}}} \nc{\CV}{{\mathcal{V}}} \nc{\CW}{{\mathcal{W}}}
\nc{\CX}{{\mathcal{X}}} \nc{\CZ}{{\mathcal{Z}}}

\nc{\LT}{{\check{T}}} \nc{\hlambda}{{\hat{\lambda}}} \nc{\hulambda}{{\widehat{-\lambda-2\rho}}}
\nc{\hLambda}{{\widehat{\lambda+\Lambda}}} \nc{\hzeta}{{\hat{\zeta}}} \nc{\gen}{{\operatorname{gen}}}
\nc{\cM}{{\check{\mathcal M}}{}} \nc{\csM}{{\check{\mathcal A}}{}} \nc{\oM}{{\overset{\circ}{\mathcal M}}{}}
\nc{\obM}{{\overset{\circ}{\mathbf M}}{}} \nc{\oCA}{{\overset{\circ}{\mathcal A}}{}}
\nc{\obA}{{\overset{\circ}{\mathbf A}}{}} \nc{\ooM}{{\overset{\circ}{M}}{}}
\nc{\osM}{{\overset{\circ}{\mathsf M}}{}} \nc{\vM}{{\overset{\bullet}{\mathcal M}}{}}
\nc{\nM}{{\underset{\bullet}{\mathcal M}}{}} \nc{\oD}{{\overset{\circ}{\mathcal D}}{}}
\nc{\obD}{{\overset{\circ}{\mathbf D}}{}} \nc{\oA}{{\overset{\circ}{\mathbb A}}{}}
\nc{\op}{{\overset{\bullet}{\mathbf p}}{}} \nc{\cp}{{\overset{\circ}{\mathbf p}}{}}
\nc{\oU}{{\overset{\bullet}{\mathcal U}}{}} \nc{\oZ}{{\overset{\circ}{\mathcal Z}}{}}
\nc{\ofZ}{{\overset{\circ}{\mathfrak Z}}{}} \nc{\ovin}{{\overset{\circ}{\on{Vin}}}}

\nc{\fa}{{\mathfrak{a}}} \nc{\fb}{{\mathfrak{b}}} \nc{\fg}{{\mathfrak{g}}}
\nc{\g}{{\mathfrak{g}}}
\nc{\fgl}{{\mathfrak{gl}}} \nc{\fh}{{\mathfrak{h}}}
\nc{\fj}{{\mathfrak{j}}} \nc{\fm}{{\mathfrak{m}}} \nc{\fri}{{\mathfrak{i}}}
\nc{\fn}{{\mathfrak{n}}} \nc{\fu}{{\mathfrak{u}}} \nc{\fp}{{\mathfrak{p}}} \nc{\frr}{{\mathfrak{r}}}
\nc{\fs}{{\mathfrak{s}}} \nc{\ft}{{\mathfrak{t}}} \nc{\fT}{{\mathfrak{T}}} \nc{\ofT}{{\overline{\mathfrak
T}}} \nc{\ofS}{{\overline{\mathfrak S}}} \nc{\fsl}{{\mathfrak{sl}}} \nc{\hsl}{{\widehat{\mathfrak{sl}}}}
\nc{\hgl}{{\widehat{\mathfrak{gl}}}} \nc{\hg}{{\widehat{\mathfrak{g}}}}
\nc{\chg}{{\widehat{\mathfrak{g}}}{}^\vee} \nc{\hn}{{\widehat{\mathfrak{n}}}}
\nc{\chn}{{\widehat{\mathfrak{n}}}{}^\vee}

\nc{\fA}{{\mathfrak{A}}} \nc{\fB}{{\mathfrak{B}}} \nc{\fC}{{\mathfrak{C}}} \nc{\fD}{{\mathfrak{D}}}
\nc{\fE}{{\mathfrak{E}}} \nc{\fF}{{\mathfrak{F}}} \nc{\fG}{{\mathfrak{G}}} \nc{\fI}{{\mathfrak{I}}}
\nc{\fK}{{\mathfrak{K}}} \nc{\fL}{{\mathfrak{L}}} \nc{\fM}{{\mathfrak{M}}} \nc{\fN}{{\mathfrak{N}}}
\nc{\frP}{{\mathfrak{P}}} \nc{\fS}{{\mathfrak S}} \nc{\fU}{{\mathfrak{U}}} \nc{\fX}{{\mathfrak{X}}}
\nc{\fY}{{\mathfrak{Y}}} \nc{\fZ}{{\mathfrak{Z}}}

\nc{\bb}{{\mathbf{b}}} \nc{\bc}{{\mathbf{c}}} \nc{\be}{{\mathbf{e}}} \nc{\bj}{{\mathbf{j}}} \nc{\bd}{{\mathbf{d}}}
\nc{\bn}{{\mathbf{n}}} \nc{\bp}{{\mathbf{p}}} \nc{\bq}{{\mathbf{q}}} \nc{\bfu}{{\mathbf{u}}}
\nc{\bv}{{\mathbf{v}}} \nc{\bx}{{\mathbf{x}}} \nc{\by}{{\mathbf{y}}} \nc{\bw}{{\mathbf{w}}}
\nc{\bA}{{\mathbf{A}}} \nc{\bB}{{\mathbf{B}}} \nc{\bC}{{\mathbf{C}}} \nc{\bF}{{\mathbf{F}}}
\nc{\bI}{{\mathbf{I}}} \nc{\bbI}{{\mathbf{I}}} \nc{\bK}{{\mathbf{K}}}
\nc{\bL}{{\mathbf{L}}} \nc{\bD}{{\mathbf{D}}}
\nc{\bH}{{\mathbf{H}}} \nc{\bM}{{\mathbf{M}}} \nc{\bN}{{\mathbf{N}}}
\nc{\bT}{{\mathbf{T}}} \nc{\bY}{{\mathbf{Y}}}
\nc{\bV}{{\mathbf{V}}} \nc{\bW}{{\mathbf{W}}} \nc{\bX}{{\mathbf{X}}}
\nc{\bP}{{\mathbf{P}}}
\nc{\bZ}{{\mathbf{Z}}} \nc{\bnu}{{\boldsymbol{\nu}}}

\nc{\sA}{{\mathsf{A}}} \nc{\sB}{{\mathsf{B}}} \nc{\sC}{{\mathsf{C}}}
\nc{\sD}{{\mathsf{D}}}
\nc{\sF}{{\mathsf{F}}} \nc{\sK}{{\mathsf{K}}} \nc{\sM}{{\mathsf{M}}}
\nc{\sO}{{\mathsf{O}}}
\nc{\sQ}{{\mathsf{Q}}} \nc{\sP}{{\mathsf{P}}} \nc{\sV}{{\mathsf{V}}}
\nc{\sW}{{\mathsf{W}}}
\nc{\sZ}{{\mathsf{Z}}} \nc{\sfp}{{\mathsf{p}}}
\nc{\sr}{{\mathsf{r}}} \nc{\sfb}{{\mathsf{b}}} \nc{\sfc}{{\mathsf{c}}}
\nc{\sd}{{\mathsf{d}}} \nc{\sk}{{\mathsf{k}}}  \nc{\sfl}{{\mathsf{l}}}

\nc{\BK}{{\bar{K}}}

\nc{\tB}{{\widetilde{\mathcal{B}}}} \nc{\tg}{{\widetilde{\mathfrak{g}}}}
\nc{\tG}{{\widetilde{G}}} \nc{\TM}{{\widetilde{\mathbb{M}}}{}} \nc{\tO}{{\widetilde{\mathbb{O}}}{}}
\nc{\tU}{{\widetilde{\mathfrak{U}}}{}} \nc{\TZ}{{\tilde{Z}}} \nc{\tx}{{\tilde{x}}} \nc{\tbv}{{\tilde{\bv}}}
\nc{\tfP}{{\widetilde{\mathfrak{P}}}{}} \nc{\tz}{{\tilde{\zeta}}} \nc{\tmu}{{\tilde{\mu}}}

\nc{\ul}{\underline} \nc{\ol}{\overline} \nc{\urho}{\underline{\rho}} \nc{\uB}{\underline{B}}
\nc{\uC}{{\underline{\mathbb{C}}}} \nc{\uc}{{\underline{c}}} \nc{\ucs}{{\underline{c},\operatorname{ss}}}
\nc{\ui}{\underline{i}} \nc{\uj}{\underline{j}} \nc{\ofP}{{\overline{\mathfrak{P}}}}
\nc{\oB}{{\overline{\mathcal{B}}}} \nc{\og}{{\overline{\mathfrak{g}}}} \nc{\oI}{{\overline{I}}}

\nc{\eps}{\varepsilon} \nc{\hrho}{{\hat{\rho}}}

\nc{\one}{{\mathbf{1}}} \nc{\two}{{\mathbf{t}}}

\nc{\Rep}{{\mathop{\operatorname{\rm Rep}}}} 
\nc{\Tot}{{\mathop{\operatorname{\rm Tot}}}} 
\nc{\Ker}{{\mathop{\operatorname{\rm Ker}}}} \nc{\Hilb}{{\mathop{\operatorname{\rm Hilb}}}}
\nc{\End}{{\mathop{\operatorname{\rm End}}}}
\nc{\Ext}{{\mathop{\operatorname{\rm Ext}}}} \nc{\Eext}{{{\mathcal{E}}xt}}
\nc{\Hom}{{\mathop{\operatorname{\rm Hom}}}} \nc{\RHom}{{\mathop{\operatorname{\rm RHom}}}}
\nc{\CHom}{{\mathop{\operatorname{{\mathcal{H}}\it om}}}} \nc{\GL}{{\mathop{\operatorname{\rm GL}}}}
\nc{\gr}{{\mathop{\operatorname{\rm gr}}}} \nc{\Id}{{\mathop{\operatorname{\rm Id}}}}
\nc{\Ind}{{\mathop{\operatorname{\rm Ind}}}} \nc{\defi}{{\mathop{\operatorname{\rm def}}}}
\nc{\length}{{\mathop{\operatorname{\rm length}}}} \nc{\supp}{{\mathop{\operatorname{\rm supp}}}}

\nc{\Cliff}{{\mathsf{Cliff}}} 
\nc{\Fl}{{\mathsf{Fl}}}
\nc{\Fib}{{\mathsf{Fib}}} \nc{\Coh}{{\operatorname{Coh}}}
\nc{\FCoh}{{\mathsf{FCoh}}}
\nc{\Mod}{{\operatorname{Mod}}}

\nc{\reg}{{\text{\rm reg}}}

\nc{\cplus}{{\mathbf{C}_+}} \nc{\cminus}{{\mathbf{C}_-}} \nc{\cthree}{{\mathbf{C}_*}} \nc{\Qbar}{{\bar{Q}}}

\nc{\bh}{{\bar{h}}} \nc{\bOmega}{{\overline{\Omega}}}

\nc{\seq}[1]{\stackrel{#1}{\sim}}

\newcommand{\beq}{\begin{equation}} \newcommand{\eeq}{\end{equation}} \renewcommand{\proof}{{\it Proof }}
\newcommand{\proofpt}{{\it Proof. }}

\newcommand{\Gm}{{\mathbb G}_m}
 
\newcommand{\TTT}{{\mathbb T}} \newcommand{\rk}{{\operatorname{rk}}}
\newcommand{\im}{{\operatorname{Im}}} \newcommand{\Aut}{{\operatorname{Aut}}}

\nc{\blambda}{{\boldsymbol{\lambda}}}
\nc{\bmu}{{\boldsymbol{\mu}}}
\newcommand\red{/\!\!/\!\!/}

\newcommand{\precc}{{\stackrel{\on{com}}{\prec}}}
\newcommand{\precg}{{\stackrel{\on{geo}}{\prec}}}
\newcommand{\preceqc}{{\stackrel{\on{com}}{\preceq}}}
\newcommand{\preceqg}{{\stackrel{\on{geo}}{\preceq}}}
\newcommand{\succeqc}{{\stackrel{\on{com}}{\succeq}}}

\renewcommand{\O}{{\mathcal{O}}}

\newcommand{\gl}{\mathfrak{gl}}

\newcommand{\la}{\lambda}

\nc{\fz}{{\mathfrak{z}}}

\nc{\Atwo}{{\mathbb{A}^2}}

\nc{\rank}{{\mathrm{rank}\,}}

\nc{\R}{{\sf R}}

\author{Roman Bezrukavnikov, Michael Finkelberg}
\title[Wreath Macdonald polynomials]{Wreath Macdonald polynomials and
the categorical McKay correspondence}

\begin{document}

\begin{abstract}
Mark Haiman has reduced Macdonald Positivity Conjecture to a statement about
geometry of the Hilbert scheme of points on the plane, and formulated a generalization
of the conjectures where the symmetric group is replaced by the wreath product
$\fS_n\ltimes (\Zet/r\Zet)^n$. He has proven the original conjecture 
by establishing  the geometric statement about
the Hilbert scheme, as a byproduct he obtained a derived equivalence
between coherent sheaves on the Hilbert scheme and coherent sheaves
on the orbifold quotient of ${\mathbb A}^{2n}$ by the symmetric group $\fS_n$.

A short proof of a similar derived equivalence for any symplectic quotient
singularity has been obtained by the first author and Kaledin
\cite{BK} via quantization in positive characteristic. In the present note we
prove various properties of these derived equivalences and then deduce
generalized Macdonald positivity for wreath products.

 \end{abstract}

\maketitle 

\centerline{\em{To the memory of Andrei Zelevinsky.}}

\tableofcontents

\section{Introduction}  The celebrated Macdonald Positivity Conjecture
asserts that the entries of the matrix expressing the transformed
Macdonald polynomials
via Schur polynomials have non-negative coefficients.

The following approach to the conjecture originated from and was completed in
the work of M.~Haiman \cite{H0}.

Let $\on{Hilb}^n(\Atwo)$ be the Hilbert scheme of $n$ points on the plane.
The action of the two-dimensional torus $T=\Gm^2$ on $\Atwo$ induces an action
of $T$ on $\on{Hilb}^n(\Atwo)$. Let $\Lambda$ be the ring of symmetric polynomials.
According to Haiman's program, one can
identify the space $\Lambda[q^{\pm1},t^{\pm1}]$ appearing in the Positivity Conjecture with
 the Grothendieck group $K^0(\Coh^T(\on{Hilb}^n(\Atwo)))$ in such a way that
   the modified Macdonald
polynomials correspond to classes of sky-scrapers at the fixed points of $T$,
 while Schur functions correspond to classes of indecomposable summands components of the {\em Procesi bundle}. The latter is a certain vector bundle on
$\on{Hilb}^n(\Atwo)$ of rank $n!$
 carrying an action of the symmetric group and a compatible action
 of the ring of polynomials in $2n$ variables.
 An explicit construction of this bundle is the crucial and laborious ingredient in Haiman's work.
 As a byproduct of this construction Haiman obtained
an  equivalence of derived categories:
\begin{equation}\label{dereq}
D^b(\Coh(\on{Hilb}^n(\Atwo))\cong D^b(\Coh^{\fS_n}({\mathbb A}^{2n}))
\end{equation}
and an  isomorphism between $\on{Hilb}^n(\Atwo)$ and the Hilbert quotient of
$\BA^{2n}$ by the action of $\fS_n$.

The existence of a rank $n!$ bundle with an action of $\fS_n$ and of the polynomial ring
which induces a  derived equivalence \eqref{dereq}
 has been proven  in another, shorter way in \cite{BK}
via quantization in positive characteristic.

The goal of the present note is to demonstrate that the information provided by
\cite{BK} is sufficient to deduce the Positivity Conjecture, bypassing the explicit
construction of the Procesi bundle.\footnote{Notice however that we do not
propose an alternative proof of
the isomorphism between $\on{Hilb}^n(\Atwo)$
 and the Hilbert quotient of $\BA^{2n}$, which
 is one of the equivalent statements
 of the $n!$ conjecture. Existence of such an isomorphism is equivalent to
the fact that each fiber
 of the rank $n!$ vector bundle is generated by an $\fS_n$ invariant vector
 under the action of the polynomial ring.
  We do not know how to deduce this property of the vector bundle by our
 methods.}
 Furthermore, Haiman~\cite{H} has suggested an extension of the
conjecture to the setting where $\fS_n$ is replaced by the wreath
product $\fS_n\ltimes (\Zet/r\Zet)^n$. We also confirm that
 previously unknown   conjecture.

We should mention that another proof of the Positivity Conjecture independent of Haiman's
construction of Procesi bundle has been obtained by I.~Gordon \cite{G_pr} based on a
result of V.~Ginzburg \cite{Gi}. This approach employs the theory of Hodge $D$-modules
to pass from modules over a noncommutative algebra to modules over its commutative degeneration; in the present article this is achieved by virtue of the $p$-center phenomenon
in the theory of $D$-modules in positive characteristic. There are other situations when 
the two apparently very different constructions lead to the same structures, see e.g. 
\cite[\S 0.2; Remark 2.2.2.(2)]{BeRi}.
 It would be very interesting to find an explanation for this parallelism.

\medskip

\subsection{}
\label{1.1}
Let us formulate our result in more detail.
We set $\Gamma=\BZ/r\BZ$.
We consider the cyclic quiver $Q$ with the set of vertices
$I=\{i_0,\ldots,i_{r-1}\}\cong \BZ/r\BZ$.
Our basic field $\sk$ will be either $\BC$ or
$\overline\BF_p$ with $p\gg0$. Given $I$-graded vector spaces $W,V$ of dimension vectors
$\lambda_0=(1,0,\ldots,0),\ \nu=(v_0,\ldots,v_{r-1})$, we consider the Nakajima
quiver varieties $\fM(W,V)\stackrel{\pi}{\longrightarrow}\fM_0(W,V)$ over $\sk$.
It is well known (see e.g.~\cite[Section~7.2.3]{H} or~\cite[Lemma~7.8]{G}) that
$\fM_0(W,V)\cong\BA^{2n}/\Gamma_n$ where $\Gamma_n$ is the wreath product
$\BZ/r\BZ\wr\fS_n$, where $n$ depends on $\nu$ as follows.
Let $\delta$ stand for
the dimension vector $(1,\ldots,1)$. Then we have a unique decomposition
$\nu=\nu_0+n\delta$ where $\nu_0$ is the content vector of a certain $r$-core
partition. Moreover, $\fM(W,V)\cong Y_{\Gamma,\nu}$, where $Y_{\Gamma,\nu}$
is a certain connected
component of the fixed point set $(\on{Hilb}^m(\BA^2))^\Gamma$ of a Hilbert
scheme of points on $\BA^2$. Here $m=v_0+\ldots+v_{r-1}$ \cite[\S 4.2]{Naka_book}.
Finally, $\pi$ is a symplectic semismall resolution of singularities.

According to the result of \cite{BK} we have an equivalence of derived
categories $D^b(\Coh( Y_{\Gamma,\nu}))\cong D^b(\Coh^{\Gamma_n}(\BA^{2n}))$
provided that characteristic of $\sk$ is zero or sufficiently large.
Let $\CE\in D^b(\Coh(Y_{\Gamma,\nu}))$ be the image of $\CO\otimes\sk[\Gamma_n]\in
\Coh^{\Gamma_n}(\BA^{2n})$ under the equivalence constructed in \cite{BK}.
Then it is shown in {\em loc. cit.} that $\CE$ is a vector bundle.
It automatically carries an action of $\Gamma_n$, thus it can canonically
be written as $\CE=\oplus \rho\otimes \CE_\rho$ where $\rho$ runs over the set of irreducible
representations of $\Gamma_n$ and $\CE_\rho$ is a vector bundle on $Y_{\Gamma,\nu}$.

Furthermore, $\CE$ carries an action of the polynomial algebra
$\sk[x_1,\dots, x_n,
y_1,\dots, y_n]=\CO(\BA^{2n})$. Let $\overline\CE_\rho$ be the image of $\CE_\rho$
in the sheaf of coinvariants $\CE/(\sum\limits_{i=1}^n x_i \CE)$,
and  $\overline\CE{}'_\rho$ be the image of $\CE_\rho$ in $\CE/(\sum\limits_{i=1}^n y_i \CE)$.

Our main result is (in the notations and conventions introduced below in
subsection~\ref{components}.)

\begin{thm}
\label{unique}
If $\rho=\rho_\blambda$ for a multi-partition $\blambda$, then the set of $T$-fixed points
in the support of $\overline\CE_\rho$ consists of
$\blambda$ and some $\bmu$ satisfying $\bmu\preceq\blambda$, while the
set of $T$-fixed points in the support of
$\overline\CE{}'_\rho$ consists of $^t\blambda$ and some $\bmu$ satisfying
$\bmu\succeq\ ^t\blambda$.
\end{thm}

The proof of the Theorem appears in subsection~\ref{210}.
It is based on the localization result of \cite{BK}, while the analysis of partial orders
on the set of fixed points follows the results of I.~Gordon~\cite{G}. Notice that Gordon
worked in a different setting, where the characteristic of the base field is zero, but
one of the deformation parameters for the deformation of $\O({\mathbb{A}}^{2n})^{\Gamma_n}$
 is also zero, so one still gets an algebra with a large center. To relate the two settings we
 include the equivalence of \cite{BK} in a family, see beginning of section \ref{new_sec}
 for further comments.

The next
statement follows from Corollary \ref{rom}.

\begin{cor}
Conjecture \cite[7.2.19]{H} holds.
\end{cor}

Notice that notations of~\cite{H} differ from ours by swapping the role of 
letters $q$ and $t$; also condition (ii) of {\em loc. cit.} appears below in
an equivalent form obtained by replacing $(-t)^{-i}$ by $(-t)^i$ and replacing
$\Lambda^i(\fh)$ by $\Lambda^{n-i}(\fh)\cong \Lambda^i(\fh^*)\otimes \Lambda^n(\fh)$.
We have $\blambda\otimes  \Lambda^n(\fh)=\ ^t \blambda$
 which accounts for appearance of  $^t\blambda$ in our Corollary \ref{rom}.
 
\subsection{Acknowledgments} We are grateful to V.~Ginzburg, I.~Gordon,
W.~van der Kallen and B.~Webster  for
some very useful discussions. Special thanks go  to I.~Losev for expert
critical comments which contributed to our understanding of the material and
influenced the exposition.

This work was conceived when the authors enjoyed the
hospitality and support of IAS at HUJI in the Fall 2010.
R.B. was partially supported by NSF grant DMS-1102434.
M.F. was partially supported by the RFBR grants 12-01-00944, 12-01-33101,
13-01-12401/13, and
the AG Laboratory HSE, RF government grant, ag. 11.G34.31.0023.


\section{Localization of modules over the invariants in the Weyl algebra}
\label{2}

\subsection{Notations}
\label{Notations}
From now on we assume that in case $\sk=\overline\BF_p$, the characteristic
$p$ is bigger than $|\Gamma_n|$. In that case, for an algebraic variety $X$
over $\sk$ we will denote by $X^{(1)}$ its Frobenius twist, and by
$\on{Fr}:\ X\to X^{(1)}$ the Frobenius morphism.
Note that in the presentation
$\fM_0(W,V)\cong\BA^{2n}/\Gamma_n$ the vector space $U=\BA^{2n}$ carries a
$\Gamma_n$-invariant symplectic form $\omega$. It gives rise to the Weyl
algebra $\sW_\sk$ acted upon by $\Gamma_n$. According to~\cite[Lemma~6.1]{BK},
for large $p$
the $\sW_\sk^{\Gamma_n}-\sW_\sk\rtimes\sk[\Gamma_n]$-bimodule $\sW_\sk$ gives rise
to the Morita equivalence $\sW_\sk^{\Gamma_n}\sim\sW_\sk\rtimes\sk[\Gamma_n]$.
We will denote $Y_{\Gamma,\nu}$ introduced in~Section~\ref{1.1} by $Y$ for short, and if we want to stress
the base field, we will denote it by $Y_\sk$.

In case $\sk=\overline\BF_p$, the Frobenius twist $Y^{(1)}$ carries a sheaf
$\sO$ of Azumaya algebras such that $H^i(Y^{(1)},\sO)=0$ for $i>0$, and
$H^0(Y^{(1)},\sO)=\sW_\sk^{\Gamma_n}$~(\cite[Lemma~6.2]{BK}, see also correction below in section \ref{erra}). From the point
of view of~\cite{EGO}, $\sO$ is the localization of the spherical subalgebra
of the symplectic reflection algebra $H_0$ of $(\Gamma_n,U)$
(rational Cherednik algebra) for the zero value of parameters.

We choose a $\Gamma_n$-invariant Lagrangian vector subspace $L\subset U$.
We will denote the categorical quotient $L/\Gamma_n\subset U/\Gamma_n$
by $\fL_0\subset\fM_0(W,V)$. We will denote the scheme-theoretic preimage
$\pi^{-1}(\fL_0)\subset\fM(W,V)=Y$ by $\fL\subset Y$. The completions of our
schemes along their subschemes will be denoted like $\widehat{Y}_\fL$.
Note that the dilation action of $\Gm$ on $\Hom(V_i,V_{i\pm1})$
descends to the $\Gm$-action
on $\fM_0(W,V)$, and then lifts to the $\Gm$-action on $Y$.

We will also need another commuting action of $\Gm$ on $Y$: the hyperbolic one, preserving
the symplectic form on $Y$. To distinguish it from the dilation action above,
we will denote this copy of $\Gm$ by ${\mathbb G}_m^h$.
So if a point $y$ of $Y$ is
represented by a quadruple $(B_i\in\Hom(V_i,V_{i+1}),\ B'_i\in\Hom(V_i,V_{i-1}),\
x_0\in V_0,\ x_0^*\in V_0^*)$, and $c\in{\mathbb G}_m^h$,
then $cy$ is represented
by $(cB_i,c^{-1}B'_i,x_0,x_0^*)$.

\subsection{Splitting bundles for  Azumaya algebras}\label{sec22}

\begin{lem}\label{zeroth}
a) There exists an Azumaya algebra $\CA$ on $({\mathbb A}^{2n})^{(1)}/\Gamma_n$,
such that the pull-back of $\CA$ to $Y^{(1)}$ is Morita equivalent to $\sO$,
while its pull-back to $({\mathbb A}^{2n})^{(1)}$ is $\Gamma_n$-equivariantly
Morita equivalent to $\sW$.

b) The Azumaya algebra $\CA$ splits on the formal neighborhood of $L^{(1)}/\Gamma_n$.
\end{lem}

\proofpt  (a) is shown in \cite[Proposition 6.5]{BK} (see also correction below in section
\ref{erra}). We will use another construction
of $\CA$, appearing  in \cite{bfgi} and recalled in the next section \ref{new_sec}. (We do not address the question of comparing
the two constructions). More precisely, we set $\sO=\bA_{1,0,..,0}$ in the notations of Proposition
\ref{exist_deform} below.

To check (b) we first check that $\CA|_{L^{(1)}/\Gamma_n}$ splits.
Since $L^{(1)}/\Gamma_n$ is smooth, in view of \cite[Corollary IV.2.6]{Mln}
(which says that the Brauer group of a regular scheme injects into the Brauer group of a dense
open subscheme) it suffices to check that the restriction of $\CA$ to
a nonempty open subset in $L^{(1)}/\Gamma_n$ splits. Let $L_0\subset L$
be the union of free $\Gamma_n$ orbits. Then splitting of $\CA|_{L^{(1)}_0/\Gamma_n}$
amounts to a $\Gamma_n$-equivariant splitting of $\sW_\sk$ on $L^{(1)}_0$.
This is provided by
the module $\CO(L)$ equipped with the natural action of $\sW_\sk$, where
elements in $L$ act by multiplication by linear functions and elements in
a fixed $\Gamma_n$-invariant Lagrangian complement $'L$ act by derivatives.
We now check  by induction that pull-back of $\CA$ to the $n$-th infinitesimal neighborhood
of $L^{(1)}/\Gamma_n$ splits; moreover, given a choice of a splitting on the $n$-th
infinitesimal neighborhood, there exists a compatible splitting on the $(n+1)$st one, which
clearly implies statement~(b).
The compatible splitting on the $(n+1)$st neighborhood exists since the obstruction lies in the second cohomology of $\CO$
on the affine variety $L^{(1)}/\Gamma_n$. $\square$

\begin{cor}
\label{first}
a) The Azumaya algebra $\sW_\sk$ splits $\Gamma_n$-equivariantly on
$\widehat{U}{}^{(1)}_{L^{(1)}}$;
$\widehat{U}{}^{(1)}_0$.

b)
 The Azumaya algebra $\sO$ splits  on
$\widehat{Y}{}^{(1)}_{\fL^{(1)}}$. $\square$

\end{cor}

Now we have the following chain of equivalences of derived categories,
cf.~\cite[Remark~6.8 and the proof of~Theorem~6.7]{BK}:
\begin{multline}
\label{mult}
 D^b\Coh(\widehat{Y}{}^{(1)}_{\fL^{(1)}})\cong
D^b(\sO|_{\widehat{Y}{}^{(1)}_{\fL^{(1)}}}-\Mod)\underset{R\Gamma}{\iso}
D^b(\sW_\sk^{\Gamma_n}|_{\widehat\fM_0(W,V){}^{(1)}_{\fL_0^{(1)}}}-\Mod)\\
\cong
D^b(\sW_\sk\rtimes\sk[\Gamma_n]|_{\widehat{U}{}^{(1)}_{L^{(1)}}}-\Mod)\cong
D^b\Coh^{\Gamma_n}(\widehat{U}{}^{(1)}_{L^{(1)}}).
%
\end{multline}
Notice that the first (respectively, the last) equivalence depends on the choice of a splitting bundle
whose existence is guaranteed by Corollary \ref{first}(b) (respectively, a).
We denote those splitting bundles by $\CE_\CL$, $\CE_L$ respectively.

\begin{prop}\label{22}
There exists an equivalence of triangulated categories
$$\Phi_\sk:D^b\Coh^{\Gmo\times(\Gmo)^h\times \Gamma_n}(U^{(1)})\iso
D^b\Coh^{\Gmo\times (\Gmo)^h}(Y^{(1)}),$$
which is compatible with the composed equivalence in \eqref{mult}
for an appropriate choice of the splitting bundles in Corollary
\ref{first}.
\end{prop}

\proofpt Let $\CA$ be as in Lemma \ref{zeroth},
and let $\CE_{res}$  be a splitting
bundle for $\sO \otimes _{\CO(Y^{(1)})} \pi^*(\CA^{op})$.
Then it is shown in \cite{BK} that the algebra $\End (\CE_{res})^{op}$ is Morita
equivalent to $\Gamma_n \# \CO({\mathbb A}^{2n})$ and
the functor $\F\mapsto R\Hom(\CE_{res},\F)$ provides an equivalence
\begin{equation}\label{eqBK}
D^b(\Coh(Y^{(1)}))\iso D^b(\End(\CE_{res})^{op}-\Mod_{fg})\cong \Coh^{\Gamma_n}(
({\mathbb A}^{2n})^{(1)}).
\end{equation}

Here the first equivalence is given by $\F\mapsto R\Hom(\CE_{res},\F)$. To describe
the second equivalence  we need to fix a $\Gamma_n$-equivariant
splitting bundle $\CE_{orb}$ for $\sW_\sk$.
Then we get:
\begin{multline*}
\End (\CE_{res})=\Gamma (\sO) \otimes _{\CO({\mathbb A}^{2n}/\Gamma_n)^{(1)} }
\Gamma(\CA)=\sW^{\Gamma_n}\otimes  _{\CO({\mathbb A}^{2n}/\Gamma_n)^{(1)}}
\Gamma(\CA)\sim \\
  \Gamma_n\# \sW   \otimes _{\CO({\mathbb A}^{2n}/\Gamma_n)^{(1)}}
\Gamma(\CA)=\Gamma_n\# \End_{\CO({\mathbb A}^{2n})^{(1)} }(\CE_{orb})\sim
\Gamma_n \#  \CO(({\mathbb A}^{2n})^{(1)}).
\end{multline*}
Here the first Morita equivalence is given by the bimodule
$\sW\otimes  _{\CO({\mathbb A}^{2n}/\Gamma_n)^{(1)}}
\Gamma(\CA)$, while the second one comes from the equivalence between
the categories of $\Gamma_n$-equivariant coherent sheaves on $({\mathbb A}^{2n})^{(1)}$
and  $\Gamma_n$-equivariant modules over $\End_{\CO({\mathbb A}^{2n})^{(1)} }(\CE_{orb})$
given by $\F\mapsto \F\otimes_{\CO} \CE_{orb}$, the latter tensor product being equipped
with the diagonal $\Gamma_n$-action.

Thus equivalence \eqref{eqBK} depends on the choice of splitting bundles
$\CE_{res}$, $\CE_{orb}$ in Corollary \ref{first}, which we assume to be fixed
from now on.  That choice determines (up to an isomorphism) a choice of splitting bundles
$\CE_L$, $\CE_{\CL}$ in the following manner.
By Lemma \ref{zeroth}(b) Azumaya algebra $\CA$ splits on the formal neighborhood
of $L^{(1)}/\Gamma_n$. The splitting bundle on that neighborhood is determined
uniquely up to twisting by a line bundle on the formal neighborhood;
since
every such bundle is easily seen to be trivial, the splitting is unique up to an isomorphism.
A diagram chase shows that under those choices the equivalences
 \eqref{mult} and \eqref{eqBK} are compatible.


We  claim that $\CE_{res}$ carries a canonical $\Gmo\times(\Gmo)^h$-equivariant structure.
This follows from (the proof of) \cite[Proposition 4.3(i)]{BK} or Lemma \ref{odin} below.



The $\Gmo\times (\Gmo)^h$-equivariant structure on $\CE_{res}$
induces a bigrading on the ring $\End(\CE_{res})$, it is easy to see
that the category of bigraded modules is canonically equivalent
to $\Coh^{\Gamma_n\times \Gmo\times (\Gmo)^h}({\mathbb A}^{2n})$.

Thus we get equivalences
\begin{equation}\label{eqbigr}
D^b(\Coh^{ \Gmo\times (\Gmo)^h}(Y^{(1)}))\cong D^b(\End(\CE_{res})^{op}-\Mod^{\mathrm{bigr}}
_{\mathrm{fg}})\cong
D^b(\Coh^{\Gamma_n\times \Gmo\times (\Gmo)^h}({\mathbb A}^{2n})^{(1)}),
\end{equation}
where $\Mod^{\mathrm{bigr}}
_{\mathrm{fg}}$ denotes the category of finitely generated bigraded modules.

It is clear from the definition that the constructed equivalence is compatible
with \eqref{eqBK}, hence with \eqref{mult}.
$\square$

 We will use the same notation $\Phi_\sk$
  for the induced compatible equivalence on categories with
 less equivariance:
 $D^b\Coh^{\Gamma_n}(\widehat{U}{}^{(1)}_{L^{(1)}}) \iso
D^b\Coh (\widehat{Y}{}^{(1)}_{\fL^{(1)}})$ and
$D^b\Coh^{\Gamma_n\times \Gmo}(\widehat{U}{}^{(1)}_{L^{(1)}}) \iso
D^b\Coh^{\Gmo} (\widehat{Y}{}^{(1)}_{\fL^{(1)}})$.

\begin{rem}
One can construct  a $\Gmo\times (\Gmo)^h$-equivariant bundle on $Y^{(1)}$
providing an equivalence $\Phi_\sk$ in a different way which is probably better suited for
generalization to other symplectic resolutions of singularities. Instead of
Lemma \ref{zeroth} one can use Corollary \ref{first}(b) to get a
tilting bundle on the formal neighborhood of $\pi^{-1}(0)$ in $Y^{(1)}$.
That bundle can be equipped with a $\Gmo\times (\Gmo)^h$-equivariant structure
by Proposition \ref{vologod} in the Appendix, which then can be shown to come
from a uniquely defined vector bundle on $Y^{(1)}$ (cf. \cite[Lemma 2.12]{BK}).
\end{rem}

\begin{lem}\label{lem_comm}
a) Assume that the splitting bundles $\CE_{res}$, $\CE_{\CL}$
and $\CE_{orb}$, $\CE_L$ are chosen compatibly as in the proof of Proposition
\ref{22}.

Then the following diagram is commutative:
$$
\begin{CD}
D^b\Coh^\Gmo_{\fL^{(1)}} (Y^{(1)})@>{\Phi_\sk^{-1}}>>
D^b\Coh^{\Gmo\times\Gamma_n}_{L^{(1)}}(U^{(1)})\\
@V{\F\mapsto R\Gamma(\F\otimes \CE_\CL)}VV @V{\F\mapsto \F\otimes \CE_L}VV \\
D^b(\sW_\sk^{\Gamma_n} -\Mod^\Gmo_{\fL_0^{(1)}})
@>{M\mapsto \sW_\sk\otimes _{\sW_\sk^{\Gamma_n} } M}>>
D^b(\sW_\sk\rtimes\sk[\Gamma_n] -\Mod^\Gmo_{L^{(1)}}),
\end{CD}
$$
where lower indices denote the set theoretic support condition on the coherent
sheaf/module/sheaf of modules.

b) For any choice of a $\Gamma_n$-equivariant splitting
$\bar{\CE}_L$ for $\sW_\sk$ on $L^{(1)}$
there exists a splitting bundle $\CE_{orb}$ compatible with $\bar{\CE}_L$.
\end{lem}

\proofpt  a) follows by diagram chase.
To check (b) notice that any two $\Gamma_n$-equivariant
splitting bundles for $\sW_\sk$ on $L^{(1)}$ differ by a twist
by a $\Gamma_n$-equivariant line bundle. It is easy to see that
 the isomorphism class of such a line bundle is determined
by the character of $\Gamma_n$ by which it acts on the fiber of the line
bundle at zero. Thus we can start with an arbitrary $\CE_{orb}$, after possibly
twisting it by a character of $\Gamma_n$ we will achieve the
 desired compatibility.
$\square$

Notice that $\sW|_{L^{(1)}}$ has a standard $\Gamma_n$-equivariant
splitting provided by the $\sW$ module $\CO(L)$ where elements in $L\subset U\cong U^*$
act by multiplication by linear functions and elements in a fixed
$\Gamma_n$-invariant
complement to $L$ act by derivatives.
By Lemma \ref{lem_comm}(b) there exists a splitting bundle $\CE_{orb}$ which is compatible
with that equivariant splitting of $\sW|_{L^{(1)}}$, from now on we assume that
the equivalence are defined using a choice of $\CE_{orb}$ satisfying that compatibility.

\medskip

Given an irreducible $\sk[\Gamma_n]$-module $\rho$, we denote by
$\CF_\rho$ the coherent sheaf $\Phi_\sk(\CO_{L^{(1)}}\otimes\rho)$.

\subsection{Lagrangian components}
\label{components}
We now fix the
  Lagrangian linear subspace $L\subset U$ given by the condition of nilpotency of
the compositions of $B_i$'s.
According to~\cite[Proposition~7.2.18]{H} or~\cite[Lemma~5.1]{G},
the ${\mathbb G}_m^h$-fixed points $Y^{{\mathbb G}_m^h}$ (recall that $Y$ was
introduced in~Section~\ref{Notations})
are naturally labeled by $r$-partitions
$\bmu=(\mu^{(1)},\ldots,\mu^{(r)})$ of total size $|\bmu|=\sum_i|\mu^{(i)}|=n$
(recall that $\nu=\nu_0+n\delta$). We will abuse notation by using the same
letter for a multi-partition and for the corresponding fixed point.

Now the Lagrangian subvariety $\fL\subset Y$ is a union of locally closed
Lagrangian subvarieties $\fL_\bmu^\circ$ labeled by the set $\CP(r,n)$ of
$r$-partitions of total
size $n$. Each piece $\fL_\bmu^\circ$ is isomorphic to the affine space $\BA^n$,
and is defined as the attracting set of the corresponding
${\mathbb G}_m^h$-fixed point. We define $\fL_\bmu$ as the closure of
$\fL_\bmu^\circ$. I.~Gordon~\cite[Section~5.4]{G} has defined a partial order
on $\CP(r,n)$ generated by the rule $\bmu\precg\blambda$ if $\bmu\ne\blambda$
and $\fL_\blambda\cap\fL_\bmu^\circ\ne\emptyset$. Moreover, he has (at least
partly) described this
order combinatorially in~\cite[Sections~6,7]{G}. In order to formulate the
description, recall
the bijection $\tau_{\mathbf s}$~\cite[6.2]{G},~\cite[7.2.17]{H} between the
set $\CP(r,n)$ of multipartitions, and the set
$\CP_{\nu_0}(|\nu_0|+rn)$ of usual partitions of size $|\nu_0|+rn$ having
$r$-core $\nu_0$. Then according to~\cite[Proposition~7.10]{G}, we have
$\bmu\precg\blambda\ \Rightarrow\ \bmu\precc\blambda\ \Leftrightarrow\
^t\tau_{\mathbf s}(^t\bmu)\lhd\
^t\tau_{\mathbf s}(^t\blambda)$, that is, one transposed partition dominates the
other. Here notation $^t\blambda$ is borrowed from \cite[\S 4.5]{G},
the multipartition $^t \blambda$ is characterized by an isomorphism 
of $\Gamma_n$ modules  $\rho_\blambda\otimes\det L\simeq\rho_{^t\!\blambda}$
(notice that the map $\blambda\mapsto \ ^t\!\blambda$ 
 is in general {\em not} an involution). 

\subsection{Localization for Verma modules for $\Gamma_n$-invariants in the Weyl algebra}
\label{LocVerm}

According to~\cite[Theorem~6.3]{BK}, the functor of global sections
$R\Gamma:\ D^b(\sO -\on{mod})\to
D^b(\sW_\sk^{\Gamma_n}-\on{mod})$ is an equivalence of categories.
We denote the quasiinverse equivalence by $M\mapsto\ ^{\on{loc}}M$ (localization).

We define {\em Verma modules}  for  $\sW_\sk\rtimes\sk [\Gamma_n]$
by: $ \sV_\bmu:=\sk[L]\otimes\rho_\bmu$ where
$\bmu$ is an $r$-multipartition of $n$, and $\rho_\bmu$ is the corresponding
irreducible $\sk[\Gamma_n]$-module. Under the Morita equivalence of
$\sW_\sk\rtimes \sk[\Gamma_n]$ and $\sW_\sk^{\Gamma_n}$, the Verma module
$\sV_\bmu$ goes to the module
$\sV_\bmu^{\Gamma_n}=\sk[L]^{\rho_\bmu}$ where the superscript
denotes the $\rho_\mu$ isotypic component. The modules $\sV_\bmu^{\Gamma_n}$ will be called the spherical Verma modules.

The key step in our argument is provided by the following:

\begin{prop}
\label{benz}
The support
$\on{supp}(\ ^{\on{loc}}\sV_{\bmu}^{\Gamma_n})$
consists of $\fL_{\bmu}^{(1)}$ and
some $\fL_\blambda^{(1)}$ such that $\blambda\precc\bmu$.
\end{prop}

The proof of the Proposition appears at the end of the next section.

\section{Spherical Cherednik algebras and localization in families}
\label{new_sec}

This section is devoted to the proof of Proposition \ref{benz}.
In order to carry
out the argument we need to include the varieties $Y$, ${\mathbb A}^{2n}/\Gamma_n$,
the ring $\sW_\sk^{\Gamma_n}$ and the sheaf
of rings $\sO$ in a family, as described in Proposition \ref{exist_deform}.

\subsection{Quiver varieties and their quantization in mixed and positive characteristic}
\label{31}
We start by presenting  the usual quiver variety construction
in a slightly  generalized setting.

We consider an arbitrary quiver $Q$ with the set of vertices $Q_0$.
We will follow the notations of~\cite[4.2]{L} for 
general Nakajima
quiver varieties. So $DQ$ is the double quiver, and given $\bv,\bd\in\BN^{Q_0}$
we consider a symplectic vector space $R(DQ,\bv,\bd)$ with a symplectic action
of $GL(\bv)$, and its Hamiltonian reduction. More precisely, given a
character $\theta:\ GL(\bv)\to\Gm$ (stability), and a central
element $\chi\in{\mathfrak z}({\mathfrak{gl}}(\bv))$
(moment level), we consider the Nakajima quiver variety
$\fM^\theta_\chi(\bv,\bd)$. It is defined in {\em loc. cit.} as a complex
quasiprojective variety, but its GIT quotient construction works over any
commutative ring
$\R$ and produces a scheme $\fM^\theta_\chi(\bv,\bd)_\R$.
If $\chi$ varies in the family ${\mathfrak z}({\mathfrak{gl}}(\bv))$,
we obtain a family of quiver varieties $\fM^\theta_{\mathfrak z}(\bv,\bd)_\R$
over ${\mathfrak z}({\mathfrak{gl}}(\bv))_\R\simeq\BA^{Q_0}_\R$.

Furthermore,
let $\sW_\sk^\hbar (\bv,\bd)$ denote the asymptotic Weyl algebra of the symplectic
vector space $R(DQ,\bv,\bd)$. Thus $\sW_\sk^\hbar(\bv,\bd)$ is  $\sk[\hbar]$-algebra
generated by $R(DQ,\bv,\bd)$ subject to the relation $[x,y]=\hbar ^2\cdot (x,y)$; it can also
be described as the Rees algebra of the Weyl algebra $\sW_\sk(\bv,\bd)$ equipped with Bernstein filtration.
The symplectic action of
$GL(\bv)$ on $R(DQ,\bv,\bd)$ gives rise to the map
$\tau:\ {\mathfrak{gl}}(\bv):=\on{Lie}GL(\bv)\to\sW_\sk^\hbar(\nu)$. For
$\chi$ such
that $\chi\cdot\bv:=\sum_{i\in Q_0}\chi_iv_i=0$ we set
$\tau_\chi(\xi):=\tau(\xi)-\hbar^2 \sum_{i\in Q_0}\chi_i\on{Tr}\xi_i$ where
$\xi=(\xi_i)_{i\in Q_0}\in{\mathfrak{gl}}(\bv)$.

The quantum Hamiltonian reduction is defined as
$\sW_\sk^\hbar(\bv,\bd)\red_\chi GL(\bv) 
:= \left( \sW_\sk^\hbar(\bv,\bd)) / \sW_\sk^\hbar(\bv,\bd)\cdot
\tau_\chi({\mathfrak{gl}}(\bv))\right)^{GL(\bv)}$ (see~\cite[4.2]{L}). Once again, this
 construction works over any ground ring $\R$ and produces a filtered $\R[\hbar]$-algebra to be denoted $A_\chi(\bv,\bd)_\R$.
This is a specialization of the
$\R[\BA^{Q_0}]$-algebra $A_{\mathfrak z}(\bv,\bd)_\R:=
 \left( \sW_\sk^\hbar(\bv,\bd)) / \sW_\sk^\hbar(\bv,\bd)\cdot
\tau({\mathfrak{sl}}(\bv))\right)^{GL(\bv)}$
at the maximal ideal  $(\chi)$ in $\R[\BA^{Q_0}]$ where
${\mathfrak{sl}}(\bv)\subset{\mathfrak{gl}}(\bv)$ is the derived subalgebra.

Finally, recall that one associates with $Q$ a Kac-Moody algebra $\fg_Q$,
and its weight $\lambda=\sum_{i\in Q_0} v_i\omega_i-\sum_{i\in Q_0}d_i\alpha_i$.
The criterion of~\cite{CB} for the moment map in the Hamiltonian reduction
construction of our quiver varieties to be flat takes especially simple
form when $\fg_Q$ is simple or affine: namely, the flatness follows if
$\lambda$ is dominant, see e.g.~\cite[2.1.4]{BL}.
In our application we only need the case of cyclic $\tilde{A}_{r-1}$-quiver
$Q$, and $\bd=(1,0,\ldots,0)$. In this case the weight $\lambda$ is dominant
iff $\bv=n\delta=(n,\ldots,n),\ n\in\BN$.

\begin{lem}\label{quivZ}
a) Consider the quiver variety $\fM^\theta_\fz(\bv,\bd)_\BZ$ over $\on{Spec}\Zet$.
There exists a finite localization $\R$ of $\Zet$ such that the base change of
$\fM^\theta_\fz(V,W)_\R$ to any algebraically closed field $\sk$
is the corresponding quiver variety $\fM^\theta_\fz(\bv,\bd)_\sk$.

b)  Assume that $\bv,\, \bd$ are such that the corresponding moment map is flat
(e.g. in case $Q$ is Dynkin or extended Dynkin the corresponding weight
$\lambda$ is dominant).
 There exists a finite localization $\R$ of $\Zet$, such that
 the associated graded of the natural filtration on $A_\fz(\bv,\bd)_\R$ is the ring of functions on the affine quiver variety $\fM_\fz^0(\bv,\bd)_\R$.
In other words, the quantum Hamiltonian reduction of
 the asymptotic Weyl algebra over $\R$ is flat over $\R[\hbar]$.

c) Assume that $\bv,\, \bd$ satisfies the assumption of (b) and $R$ satisfies both (a) and (b).
Then we have
$$\on{gr}(A_\fz(\bv,\bd)_\sk)=\O(\fM_\fz^0(\bv,\bd)_\sk ), $$
$$\on{gr}(A_\chi(\bv,\bd)_\sk)=\O(\fM_\chi^0(\bv,\bd)_\sk ), $$
 for any field $\sk$ with a homomorphism from $R$ and $\chi\in \fz_\sk^*$.

 \end{lem}

\proofpt a)
Taking the preimage of $\fz$ under the moment map commutes with any base change.
After base change to a finite localization of $\Zet$
the categorical quotient by the action of $GL(\bv)$ exists and commutes with any base change
by~\cite[Theorem~33]{FvdK}.

b)  The quantum Hamiltonian reduction is defined as the module  of
$GL(\bv,\R)$ invariants in a quotient of the Weyl algebra by a left ideal.
Thus it carries a natural filtration whose associated graded is a submodule in the module
of $GL(\bv,\R)$ invariants of the associated graded of the quotient by the left ideal.
As pointed out above, the condition on the weight $\la$ implies that the moment map is flat.
Flatness of the moment map implies in turn that taking the quotient by the ideal commutes with taking associated graded, thus we see that the associated graded of the filtration on the quantum
Hamiltonian reduction is a subring in the classical Hamiltonian reduction.
 Furthermore,
the ambient ring here is finitely generated and the embedding becomes an
isomorphism upon base change to $\Ce$.
 It follows that the two rings become equal upon base change to a finite localization of $\R$.

c) follows since there is a natural injective
 map from the left hand side to the right hand side, statements (a,b) together show these maps are surjective.
\qed

In the rest of this subsection we work over a fixed field $\sk=\overline\sk$ of characteristic $p>0$.

Recall that for an algebraic variety $X$ over $\sk$
we denote by $X^{(1)}$ its Frobenius
twist and  $\on{Fr}:\ X\to X^{(1)}$ is the Frobenius morphism.

\begin{lem}\label{char_p_quiv}  Assume that $p>\bv_i$ for all $i$.

a) The ring $A_{\mathfrak z}(\bv,\bd)_\sk$   
has a natural structure of a ring over
 $\O(\fM^0_{\mathfrak z}(\bv,\bd)^{(1)}
 \times_{(\BA^{Q_0})^{(1)}}(\BA^1_\hbar\times\BA^{Q_0}))$, where
we used the map
$AS_Q:\ (\hbar, \chi_0,\dots, \chi_{r-1})\mapsto (\chi_0^p-\hbar^{p-1}\chi_0,
\dots,\chi_{r-1}^p-\hbar^{p-1}\chi_{r-1}),\ \BA^1_\hbar\times\BA^{Q_0}\to
(\BA^{Q_0})^{(1)}$.

b) Reducing the structure described in (a) at $\hbar =0$
we get a map $\O(\fM^0_{\mathfrak z}(\bv,\bd)^{(1)}\to\O(\fM^0_{\mathfrak z}(\bv,\bd))$;
this map equals pull-back under the Frobenius morphism.
\end{lem}

\proofpt
To check (a)  we address following~\cite[Section~3]{bfgi} the special features of
Hamiltonian reduction in positive characteristic.

For a connected smooth linear algebraic group $G$ over $\sk$ we have an exact
sequence of groups $1\to G_1\to G\stackrel{\on{Fr}}{\longrightarrow}G^{(1)}\to1$
where $G_1$ stands for the Frobenius kernel. The Lie algebra $\fg$ of $G$
is equipped with a natural structure of a $p$-Lie algebra, and its universal
enveloping algebra $U(\fg)$ contains the $p$-center $\fZ(\fg)$.

We denote by ${\mathbb X}^*(G)$ the lattice of characters  $G\to \Gm$ and
we set ${\mathbb X}^*(\fg)={\mathbb X}^*(G) \otimes \sk$.
We have  ${\mathbb X}^*(\fg)\subset  \fg^*$.


For $\chi \in X^*(\fg)$ we denote by $I_{\chi}$ the kernel of the corresponding homomorphism
$U(\fg)\to\sk$. We set
$I_{\chi}^{(1)}:=I_{\chi}\cap\fZ(\fg)$, a maximal ideal of $\fZ(\fg)$.
We let
 $\fu_{\chi}(\fg)$ denote the quotient of $U(\fg)$ by the two-sided ideal
generated by
$I^{(1)}_{\chi}$. The image of $I_{\chi}$ in $\fu_{\chi}(\fg)$
is denoted by $\fri_{\chi}\subset\fu(\fg)$.

We now  take $G=GL(\bv)$ and we denote its
Lie algebra by ${\mathfrak{gl}}(\bv)$. Thus $X^*(\fg)=\fz^*$
(where we used the assumption $p>\bv_i$).

 The center $\fZ(\bv,\bd)$ of the Weyl
algebra $W_\sk(\bv,\bd)$ equals $\sk[R(DQ,\bv,\bd)^{(1)}]$, while that of
 $W_\sk^\hbar(\bv,\bd)$ equals $\sk[\BA^1_\hbar\times R(DQ,\bv,\bd)^{(1)}]$.

Thus  $W_\sk^\hbar(\bv,\bd)$ can be viewed as a coherent sheaf of algebras
over $\BA^1_\hbar\times R(DQ,\bv,\bd)^{(1)}$.

The quantized moment map is the homomorphism ${\mathfrak{gl}}(\bv)\to W_\sk^\hbar (\bv,\bd)$.
Its restriction to $\fZ({\mathfrak{gl}}(\bv))$ corresponds to the map of
spectra $\mu_\hbar^{AS}: \BA^1_\hbar \times R(DQ,\bv,\bd)^{(1)}\to
({\mathfrak{gl}}(\bv)^*)^{(1)}$ given by
$\mu_\hbar^{AS}(\hbar, x)=AS_Q(\hbar,\mu(x))$,
where $\mu:\ R(DQ,\bv,\bd) \to{\mathfrak{gl}}(\bv)^*$ is the moment
map~\cite[Proposition~$2.1'$]{T}.

We obtained a map from $\O
(\mu_\hbar^{AS})^{-1}(X^*(\fg))/\!/GL(\bv)^{(1)})=
\O(\fM^0_{\mathfrak z}(\bv,\bd)^{(1)}
 \times \BA^1_\hbar)$ to the center of the quantum Hamiltionian reduction $\fZ\left(\left(\sW_\sk^\hbar(\bv,\bd)/\on{Im}({\mathfrak{sl}}(\bv)\right)^{GL(\bv)}\right)$. We also have a map from
 $\sk[\fz]$ to $\fZ\left(\left(\sW_\sk^\hbar(\bv,\bd)/\on{Im}({\mathfrak{sl}}(\bv)\right)^{GL(\bv)}\right)$.

 Combining the two maps together we get
 a map
 \begin{equation}
\label{map} \O(\fM^0_{\mathfrak z}(\bv,\bd)^{(1)}
 \times_{(\BA^{Q_0})^{(1)}}(\BA^1_\hbar\times\BA^{Q_0}))
 \to \fZ\left(\left(\sW_\sk^\hbar(\bv,\bd)/\on{Im}({\mathfrak{sl}}(\bv)\right)^{GL(\bv)}\right),
 \end{equation}
 where the fiber product involves the map
 $\mu_\hbar ^{AS}$, this follows from the definition of $\mu_\hbar^{AS}$. This proves (a).


The specialization of the $p$-center of the Weyl algebra at $\hbar=0$ produces
the subring $\O(R(DQ,\bv,\bd))^{(1)} \overset{\on{Fr}^*}{\longrightarrow}
\O(R(DQ,\bv,\bd))$. Applying Hamiltonian reduction we get (b).
\qed

\begin{rem} It is not hard to see that at least for large $p$ and
$\bv,\, \bd$ such that the corresponding moment map is flat,
the map from $\O(\fM^0_{\mathfrak z}(\bv,\bd)^{(1)}
 \times_{(\BA^{Q_0})^{(1)}}(\BA^1_\hbar\times\BA^{Q_0}))$ to the center
 of $A_{\mathfrak z}(\bv,\bd)_\sk$ is an isomorphism.
\end{rem}

We now combine  Hamiltionian reduction in positive characteristic as recalled
in the proof of~Lemma~\ref{char_p_quiv} with the GIT quotient procedure.

More precisely, we consider  $\tilde \CA:=(\sW^\hbar_\sk(\bv,\bd)^{(1)}/
\sW_\hbar^\sk(\bv,\bd)^{(1)}\cdot{\mathfrak{sl}}(\bv))^{GL(\bv)_1}$.
This is a sheaf of associative rings with a $GL(\bv)^{(1)}$ action.
As in the proof of the previous Lemma we have maps
from $\O(\mu^{-1}(X^*(\gl(\bv)))^{(1)}$ and from $\sk[\fz\oplus \sk\cdot \hbar]$
to the center $\fZ(\tilde \CA)$; together they yield a  $GL(\bv)^{(1)}$-equivariant  map
 \begin{equation}\label{to_center}
 \O(\mu^{-1}(X^*(\gl(\bv)))^{(1)}\times_{(\BA^{Q_0})^{(1)}}
 (\BA^1_\hbar\times\BA^{Q_0}))
 \to \fZ(\tilde \CA).
\end{equation}
Thus we can view $\tilde \CA$ as a $GL(\bv)^{(1)}$-equivariant sheaf over
 $\mu^{-1}(X^*(\gl(\bv)))^{(1)}\times_{\BA^{Q_0}}
 (\BA^1_\hbar\times\BA^{Q_0}))$.

We now fix a generic (i.e. lying off the walls,
see~\cite[2.3]{NQ}) stability parameter $\theta$
and let $\mu^{-1}(X^*(\gl(\bv)))_\theta\subset \mu^{-1}(X^*(\gl(\bv)))$ be the
 subset of $\theta$-stable points.
The quotient $\mu^{-1}(X^*(\gl(\bv)))_\theta/GL(\bv)$ is by definition the
family of quiver varieties $\fM_\fz^\theta$.
We set $(\fM_\fz^\theta)^{AS}=
(\mu^{-1}(X^*(\gl(\bv)))_\theta^{(1)}\times_{\BA^{Q_0}}
(\BA^1_\hbar\times\BA^{Q_0})))/GL(\bv)$,
clearly $(\fM_\fz^\theta)^{AS}=\fM_\fz^\theta\times_{\BA^{Q_0}}
 (\BA^1_\hbar\times\BA^{Q_0}))$.

\begin{lem}\label{char_p_quiv_thet}
a) The sheaf $\tilde \CA$ descends to
a sheaf of rings $\CA_\fz^\theta$ over
$(\fM^\theta_\fz)^{AS}$.
The restriction of this sheaf to the open part $\hbar\ne 0$ is an Azumaya
algebra, while its restriction
to the closed subvariety $\hbar=0$ is isomorphic to $\on{Fr}_*(\O)$.

b) Assume that $(\bv,\bd)$ satisfy the condition of Lemma \ref{quivZ}(b)
 and $p\gg 0$. Then $R\Gamma(\CA_\fz^\theta)=A_\fz(\bv,\bd)$,
 while for every $\chi$ we have
$R\Gamma(\CA_\chi^\theta)=A_\chi(\bv,\bd)$.
\end{lem}

\proofpt
a) The action of $GL(\bv)$ on $\mu^{-1}(X^*(\gl(\bv)))_\theta$
 is free. This is well known in characteristic zero
(see e.g.~\cite[Proposition~2.6]{Y}),
and the case of positive characteristic is checked by a similar
argument.\footnote{In {\em loc. cit.} it is checked that the action
of $GL(\bv)$ on $\mu^{-1}(X^*(\gl(\bv)))_\theta$ has trivial stabilizers.
In general this does not imply that the action is free, see \cite[Example 0.4]{Mum}
(recall that an action of $G$ on $X$
is free if the map $(a\times pr_2):G\times X\to X\times X$
is a closed embedding, where $a$ is the action and $pr_2$ is the second projection).
However, here we are dealing with an action of $GL(\bv)=\prod GL(\bv_i)$
on a locally closed subvariety in the linear representation $R(DQ,\bv,\bd)$.
Given $x,\, y\in R(DQ,\bv,\bd)$ and $g=(g_i)\in GL(\bv)$ the equation $g(x)=y$ is linear
in the matrix coefficients of $g_i\in GL(\bv_i)$. By standard linear algebra, locally in $x$
we have an algebraic expression for these matrix coefficients in terms of coordinates of $x$, $y$
provided that a unique  solution exists. This shows that the action is free if the stabilizers
are trivial.}

 For a free action the factorization is a principal $G$-bundle, so descent
 theory applies, see~\cite[Proposition 0.9]{Mum}. Applying this
 to $\tilde \CA$ 
  we get a sheaf of rings on $(\fM^\theta_\fz)^{AS}$. The first property is established
 as in \cite{bfgi} and the second one is clear from the construction.

  Under the conditions of (b) the higher direct image of the structure
sheaf under the maps $\fM_\fz^\theta\to \fM_\fz^0$,
$\fM_\chi^\theta\to \fM_\chi^0$
 vanish:  in characteristic zero,
this follows from Grauert-Riemmenschneider Theorem, the case of large positive
characteristic follows since the support of the higher direct image of $\O$ under  the morphism
 of (the family of) quiver varieties over $\Zet$ is a $\Gm$-invariant closed subset which does not
 intersect the fiber over the generic point of $\on{Spec}(\Zet)$, thus it is contained in the preimage
 of a finite subject in $\on{Spec}(\Zet)$.

Hence the higher direct image of $\on{Fr}_*(\O_{\fM_\fz^\theta})$,
$\on{Fr}_*(\O_{\fM_\chi^\theta})$  also vanish, thus the same is true for $\CA$.
Thus its global sections is a {\em Cohen-Macaulay} $\O_{\fM_\fz^0}$-
(resp. $\O_{\fM_\chi^0}$-) module.
On the open part $U$ where
the map $\fM^\theta \to \fM^0$ is an isomorphism it is isomorphic to the restriction of $A_\fz(\bv,\bd)$
(respectively, $A_\chi(\bv,\bd)$.
This latter algebra is also a Cohen-Macaulay module in view of Lemma \ref{quivZ}(b).
The complement to $U$ has codimension at least two; more precisely,
$\fM^\theta_\chi \to \fM^0_\chi$ is semismall, while  $\fM^\theta_\fz \to \fM^0_\fz$
is small (this is known in characteristic zero~\cite{NQ},
hence in large positive characteristic).
Thus we get (b).
 \qed

\subsection{The family of Cherednik algebras and their localization}
\label{loca}
%
Recall that  the algebras $\sW_\sk^{\Gamma_n}\subset\sW_\sk\rtimes\sk[\Gamma_n]$
fit into the $(r+1)$-parametric family of (spherical) rational Cherednik algebras
$e\bH e\subset\bH$, see e.g.~\cite[6.1]{L}. The parameters are denoted by
$\hbar,c_0,c_1,\ldots,c_{r-1}$ (in {\em loc. cit.} $c_0$ is denoted $k$).
The case of the Weyl algebra corresponds to
$\hbar=1,c_0=c_1=\ldots=c_{r-1}=0$.

Recall  also that a spherical rational Cherednik algebra specialized at
$\hbar=0$ is commutative, thus we get a
family of affine {\em cyclic Calogero-Moser spaces} $\bZ$,  over the space
$\BA^r_c$ with coordinates $c_0,\dots, c_{r-1}$, namely $\bZ=\on{Spec}(e\bH e)$,
and the fiber $\bZ_{\underline{c}}=\on{Spec}(e\bH_{(0,\underline{c})} e)$.

\bigskip

We now present a description of these algebras as a quantization of an affine quiver variety
constructed via quantum Hamiltonian reduction.

In order to be able to apply Lemma \ref{quivZ}
we use   quiver data different from the one in \S \ref{1.1}.
Namely, the smooth quasiprojective Nakajima quiver variety $\fM(V,W)$ is defined as
a GIT quotient with respect to the stability vector $\theta=(1,\ldots,1)$,
that is the character of $G_\nu=GL(\nu):=\prod_{i\in I}GL(V_i)$ equal to
$\theta(\underline{g})=\prod_{i\in I}\det(g_i)^{-1}$. We choose an element
$\sigma$ of the $\tilde{A}_{r-1}$ affine Weyl group such that
$\sigma(\lambda_0-n\delta)=\lambda_0-\nu$, and consider
$\theta_\nu=\sigma^{-1}(\theta)$. Then the works of Lusztig, Nakajima and Maffei
on reflection functors provide an isomorphism $\fM(V,W)=\fM^\theta(V,W)\iso
\fM^{\theta_\nu}(V',W)$ where $\underline{\dim}V'=n\delta$, and
$\ul{\dim}W=\lambda_0=(1,0,\ldots,0)$, see e.g. \cite[\S 2.1.3]{BL} and references therein.
These references treat the case of a coefficient field of characteristic zero, but the same argument
applies for a field of large positive characteristic.

 From now on we will
replace the stability condition $\theta$ by $\theta_\nu$, and $V_\bullet$ by
$V'_\bullet$. Thus we have:
\begin{equation}\label{Yis}
Y=\fM^{\theta_\nu}(V',W).
\end{equation}



\begin{lem}\label{qHam}
For an appropriate choice of isomorphism
$\BA^r_c\to\BA^I,\ \underline{c}\to\chi(\underline{c})$, we have
$e\bH_{(\hbar,\ul{c})}e\cong\sW_\sk^\hbar ({n\delta})\red_{\chi(\ul{c})}G_{n\delta}$.

More precisely,  set $\eta=\exp(2\pi i/r)\in\BZ[\sqrt[r]{1}]\subset\BC$.
Choose a prime ideal ${\mathfrak p}\subset\BZ[\sqrt[r]{1}]$ over $p$,
and keep the same notation for the reduction of $\eta$ in
$\BZ[\sqrt[r]{1}]/{\mathfrak p}\subset\overline\BF_p$.

Then $\chi(\ul{c})_0:=c_0+r^{-1}(1-r-\sum_{m=1}^{r-1}c_m),\
\chi(\ul{c})_l:=r^{-1}(1-\sum_{m=1}^{r-1}\eta^{ml}c_m)$ for $1\leq l\leq r-1$.

\end{lem}
\proofpt 
A similar statement in characteristic zero is shown in~\cite[Theorem~3.13]{G1}.
To deduce the case of positive characteristic we use the following easy general statement.
Let $\R$ be a commutative ring finitely generated over $\Zet$.
Assume that two finitely presented algebras 
over $\R$
have isomorphic base changes to $\Ce$; then they become isomorphic
after base change to some finitely generated commutative flat $\R$-algebra, in particular
their base change to an algebraically closed field of almost any prime characteristic. The two algebras in question are the spherical rational DAHA and the algebra obtained by
 quantum Hamiltonian reduction from the Weyl algebra over $\Zet$.

To make sure they are finitely presented it suffices to check that
they (or rather their base change to a finite localization of $\R$)
 admit a filtration with a commutative finitely generated associated graded.
For the spherical Cherednik algebra this is clear, as it has a filtration whose
associated graded is the ring of $\Gamma_n$-invariants in the symmetric algebra.
For the quantum Hamiltonian reduction this follows from Lemma \ref{quivZ}.
\qed


We now assume that $char(\sk)= p\gg 0$. We let $\bZ^{AS}=\bZ_\sk^{(1)}\times_{(\BA^r)^{(1)}}
(\BA^1_\hbar\times\BA^r)$, where the map
$AS_Q:\BA^1_\hbar\times\BA^r\to (\BA^r)^{(1)}$ is used.

\begin{prop} \label{exist_deform}
Assume that $char(\sk)= p\gg 0$.
 There exists
 a smooth family $\bY$ of algebraic varieties over the base $\BA^{r+1}=\BA^1_\hbar\times\BA^r_c$ with
coordinates $\hbar,c_0,c_1,\ldots,c_{r-1}$,
and a flat sheaf of algebras $\bA$ 
 on $\bY^{(1)}$ with the following properties.

\begin{enumerate}
\item The fiber $\bY_{(1,0,\ldots,0)}$ equals $Y$ (introduced in~Section~\ref{Notations});

\item We have a proper map from $\bY$ to $\bZ^{AS}$.
This map is an isomorphism over an open subset $(\bZ^0)^{AS}\subset \bZ^{AS}$,
which contains the fiber $\bZ_{\ul{c},\sk}$ for a  generic $\underline{c}$.



\item We have $\bA_{(0,\underline{c})} 
\cong \on{Fr}_*\CO_{\bY_{(0,\underline{c})}}$, while the restriction of $\bA$ to the open subspace
$\hbar\ne 0$ is an Azumaya algebra.

\item $\Gamma(\bY^{(1)},\bA)=
e\bH e$.
\item The dilation action of $\Gm$ on parameters, $t:(\hbar, c_0,\dots, c_{r-1})
\mapsto (t\hbar, tc_0,\dots, tc_{r-1})$ lifts to an action on $\bY$; the sheaf $\bA$
admits a natural $\Gm$-equivariant structure.
\end{enumerate}
\end{prop}

\proofpt
We let $\bY=(\fM^{\theta_\nu}_\fz(n\delta,\lambda_0)^{AS}$, and we consider the sheaf of algebras
$\bA=\CA^{\theta_\nu}_\fz$.

Then part (a) follows from \eqref{Yis}, part (d) from Lemma \ref{char_p_quiv_thet}(b)
and Lemma \ref{qHam}, part (c) from Lemma \ref{char_p_quiv_thet}(a), and (e) is clear from the construction.
Part (b) follows from part (d) and Lemma
\ref{qHam}. 
\qed


\subsection{Verma modules}
\label{irr O}
According to~\cite[Theorem~7.2.4]{bfgi}, the functor of global sections
$R\Gamma:\ D^b(\bA_{\hbar, \chi(\ul{c})}-\on{mod})\to
D^b(e\bH_{(\hbar,\ul{c})}e-\on{mod})$ is an equivalence of categories
provided that $\bA_{\hbar, \ul{c}}$ has finite homological
dimension; parameters $(\hbar, \ul{c})$ satisfying this property will be called good.

Good parameters form a dense Zariski open subset in the space of parameters.
While it is possible to describe this set explicitly (for large $p$) based on results
of \cite{DuGr},  we will only need the following easily available information:
$(\hbar=1, \ul{c}=0)$ is good; for a generic $\ul{c}$ the pair $(\hbar,\ul{c})$ is good for any $\hbar$.

We denote the quasiinverse equivalence by $M\mapsto\ ^{\on{loc}}M$ (localization).

We will apply this notation both for a specific good parameter and for a family of such.


Recall  Verma modules for $\sW_\sk\rtimes\sk[\Gamma_n]$ and $\sW_k^{\Gamma_n}$
introduced in \ref{LocVerm}. The same construction applied to the Cherednik
algebra with variable parameters produces the Verma $\bH$-modules
$\bV_\bmu$, and the spherical Verma $e\bH e$-modules $\bV_\bmu^{\Gamma_n}$.
The specialization of $\bV_\bmu^{\Gamma_n}$ to
$(\hbar,\ul{c})\in\BA^1_\hbar\times\BA^r_c$ will be denoted by
$\bV_{\bmu,(\hbar,\ul{c})}^{\Gamma_n}$: a spherical Verma module over
$e\bH_{(\hbar,\ul{c})}e$. So for example
$\bV_{\bmu,(1,0,\ldots,0)}^{\Gamma_n}=\sV_\bmu^{\Gamma_n}$.


\begin{lem}
\label{exact} We have
$^{\on{loc}}\sV_{\bmu}^{\Gamma_n}\in\bA_{\chi(\ul{0})}-\on{mod}\subset
D^b(\bA_{\chi(\ul{0})}-\on{mod})$, i.e. $^{\on{loc}}\sV_\bmu^{\Gamma_n}$ is an
$\bA_{\chi(\ul{0})}$-module (as opposed to a complex of modules).
\end{lem}

\proofpt
We choose a complementary $\Gamma_n$-invariant Lagrangian subspace $'L:\
'L\oplus L=U$. The composed projection $Y\to\BA^{2n}/\Gamma_n\to\ 'L/\Gamma_n=\BA^n/\Gamma_n$ is flat since all the fibers are of the same
 dimension $n$. Hence a
sequence $e_1,\ldots,e_n$ of generators of $\sk[\BA^n/\Gamma_n]$ is a
regular sequence for $\CO_{Y}$, hence so is
 the sequence
$e_1^p,\ldots,e_n^p$.
Thus $e_1^{(1)}, \ldots,e_n^{(1)}$
 is a regular sequence for $\on{Fr}_*\CO_{Y}$ on $Y^{(1)}$,
 where $e_i^{(1)}$ is the function $e_i^p$ viewed as a function on $Y^{(1)}$.
  Using Proposition \ref{exist_deform}(c,e) we see that $e_1^{(1)}, \ldots,e_n^{(1)}$ is a regular sequence for the Azumaya algebra
$\sO=\bA_{\chi(\ul{0})}$ on $Y^{(1)}$. According to Corollary~\ref{first},
on the formal neighborhood of the central fiber, $\bA_{\chi(\ul{0})}\simeq
\End(E)$ splits. Since $E$ contains $\CO_{Y^{(1)}}$ as a direct
summand, $\End(E)$ contains $E$ as a direct summand, hence  $e_1^{(1)}, \ldots,e_n^{(1)}$ forms a regular sequence for $E$ as well.

To finish the proof observe that $\sV_\bmu^{\Gamma_n}$ is direct summand in $\sW^{\Gamma_n}\otimes
 _{\sk[e_1,\dots, e_n]}\sk$: in characteristic zero this is clear since
category $\CO$ for $c=0$ is semi-simple, while the Morita equivalence
$\sW^{\Gamma_n}\sim \sW$ sends $\sW^{\Gamma_n}\otimes
 _{\sk[e_1,\dots, e_n]}\sk$ to $\sW\otimes
 _{\sk[e_1,\dots, e_n]}\sk= \sW\otimes _{\sk['L]}(\sk['L]/(e_1,...e_n))$
 which carries a filtration whose associated graded is a direct sum of Verma modules, each
 one appearing with a nonzero multiplicity. It follows that
  $\sV_\bmu^{\Gamma_n}$ is direct summand in $\sW^{\Gamma_n}\otimes
 _{\sk[e_1,\dots, e_n]}\sk$ also when $\sk$ has a large positive characteristic.
 Now it is clear that
  $\sW^{\Gamma_n}\otimes
 _{\sk[e_1,\dots, e_n]}\sk=R\Gamma(\sO\Lotimes _{\sk[e_1,\dots, e_n]}\sk)$ thus $^{\on{loc}}\sV_{\bmu}^{\Gamma_n}$
is a direct summand in $\sO\Lotimes _{\sk[e_1,\dots, e_n]}\sk$.
\qed

\subsection{Supports of Verma modules}
\label{supp}
We  need some properties of the bijection between
irreducible $\Gamma_n$ modules and $(\Gm)^2$-fixed points, which will be
presently deduced
from the results of~\cite{G0}. In fact, {\em loc. cit.} establishes similar
results for some spherical
Cherednik algebras over a field of zero characteristic.
Passage between the present setting
and that of {\em loc. cit.} is done in several steps,
which schematically can be depicted as follows:

 \Red
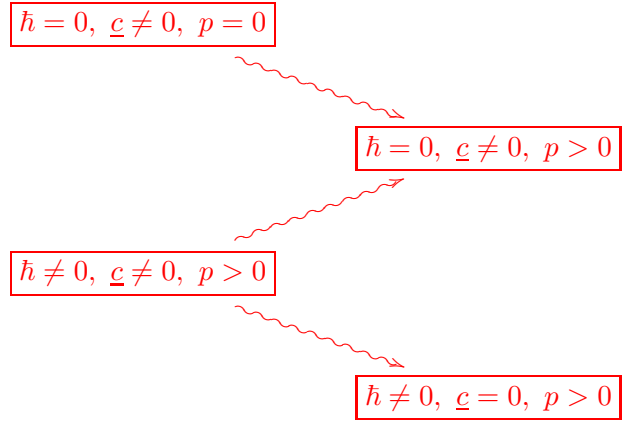
{
$$\xymatrix{
\framebox{$\hbar=0,\ \ul{c}\ne0,\ p=0$} \ar@{~>}[dr] \\
& \framebox{$\hbar=0,\ \ul{c}\ne0,\ p>0$} \\
\framebox{$\hbar\ne0,\ \ul{c}\ne0,\ p>0$} \ar@{~>}[ur] \ar@{~>}[dr]\\
& \framebox{$\hbar\ne0,\ \ul{c}=0,\ p>0$}
}$$
}

Here the first two degenerations are analyzed in Lemma \ref{benzrou}, while the last
one will be treated using Lemmas \ref{exact}, \ref{loseu}.


\bigskip

We start with a geometric Lemma on the behavior of our Lagrangian subvarieties under the
degeneration $\ul{c}\to 0$; it will be applied for varieties over a field of characteristic $p\gg 0$.

The hyperbolic ${\mathbb G}_m^h$ of
Section~\ref{components}
acts on $\bY$, and the ${\mathbb G}_m^h$-fixed points in a general fiber
are uniformly numbered by the set of $r$-partitions of $n$, according
to~\cite[3.8--3.10]{G}. As in Section~\ref{components}, we define
$\fL^\circ_{\blambda,\underline{c}}\subset\bZ_{\underline{c}}$
as the attracting set of the corresponding ${\mathbb G}_m^h$-fixed point.
Contrary to the situation of Section~\ref{components}, for a general
$(\underline{c})$ the Lagrangian subvariety
$\fL^\circ_{\blambda,\underline{c}}\subset\bZ_{\underline{c}}$
is closed, i.e. equals its closure $\fL_{\blambda,\underline{c}}$.

As $\underline{c}$ goes to $\underline{0}$, the Lagrangian subvariety
$\fL_{\bmu,\underline{c}}$ degenerates into the (closed) Lagrangian
subvariety contained in $\fL_{\bmu}\cup\bigcup_{\blambda\precc\bmu}\fL_\blambda$.
In effect, let $\bL_\bmu\subset\bY$ be the closed irreducible subvariety of
$\bY$ whose fiber over general $(\hbar,\ul{c})$ is
$\fL_{\bmu,\ul{c}^p-\hbar^{p-1}\ul{c}}$. The following lemma is well known to the
experts.\footnote{We thank Ben Webster who explained it to us.}

\begin{lem}
\label{loseu}
The intersection $\bL_\bmu\cap Y$ contains
$\fL_{\bmu}\cup\bigcup_{\blambda\precg\bmu}\fL_\blambda$ and is contained in
$\fL_{\bmu}\cup\bigcup_{\blambda\precc\bmu}\fL_\blambda$.
\end{lem}

\proofpt
The intersection $\bL_\bmu\cap Y$ is Lagrangian in $Y$ since it is
contained in $\pi^{-1}\fL_0=\fL$.
Since $\bL_\bmu\cap Y$ is closed and contains $\fL_\bmu$, it must contain
$\fL_{\bmu}\cup\bigcup_{\blambda\precg\bmu}\fL_\blambda$. Conversely, if
$\fL_\blambda\subset\bL_\bmu\cap Y$, then $\blambda\precc\bmu$. To see this
let us consider first the case $r=1$ (no cyclic group). Then both the
geometric and combinatorial orders are nothing but the dominance order on
partitions.
The parameter $\underline{c}$ is just one number $c$, and $\bY$ is just the
family degenerating the Calogero-Moser space into the Hilbert scheme.
The equivariant Borel-Moore \'etale homology groups
$H^{\BG_m^h}_{BM}(\bY_{(\hbar,c)})$
are canonically identified for all the fibers of $\bY$. The fundamental
classes of $\fL_{\mu,c^p-\hbar^{p-1}c}$ form a basis in $H^{\BG_m^h}_{BM}(\bY_{(\hbar,c)})$.
For $c^p-\hbar^{p-1}c\ne0$, when $\fL_{\mu,c^p-\hbar^{p-1}c}$ is closed, its
fundamental class in the localized equivariant Borel-Moore homology is
proportional to the class of the fixed point $\mu$. On the other hand,
under the above identification, this class equals a class in
$H^{\BG_m^h}_{BM}(Y)$ supported in $\bL_\mu\cap Y$ (and proportional to the
fixed point class). Now under the well known
identification of $H^{\BG_m^h}_{BM}(Y)$ with the degree $n$ part of the Fock
space of symmetric functions, the fixed point classes correspond to the
Schur functions, while the fundamental classes of $\fL_\lambda$ correspond
to the monomial functions (see e.g.~\cite{nak}). The transition matrix between
these two bases is upper triangular in the dominance order, hence we are done.

For general $r$, recall that $Y=Y_{\Gamma,\nu}$ is a connected component of the
fixed point set $(\on{Hilb}^m(\BA^2))^\Gamma$ where $m=v_0+\ldots+v_{r-1}$.
A fixed point $\blambda\in Y$ goes to the fixed point
$^t\tau_{\mathbf s}(\blambda)$ (a partition of $m$; notations of~\ref{components})
of $\on{Hilb}^m(\BA^2)$. If $\blambda$ lies in $\bL_\bmu\cap Y$, then
$^t\tau_{\mathbf s}(\blambda)$ lies in $\bL_{^t\tau_{\mathbf s}(\bmu)}\cap
\on{Hilb}^m(\BA^2)$, hence $^t\tau_{\mathbf s}(\blambda)\lhd\
{}^t\tau_{\mathbf s}(\bmu)$, that is $\blambda\preceqc\bmu$. \qed

\begin{lem}
\label{benzrou}
For a general $(\underline{c})$, the support
$\on{supp}(\ ^{\on{loc}}\bV_{\bmu,(\hbar,\underline{c})}^{\Gamma_n})$ equals
$\fL_{\bmu,(\underline{c}^p-\hbar^{p-1}\underline{c})}^{(1)}$.
\end{lem}

\proofpt As the Verma module $\bV_{\bmu,(\hbar,\underline{c})}$ as well as its
spherical counterpart $\bV_{\bmu,(\hbar,\underline{c})}^{\Gamma_n}$
is indecomposable,
the support $\on{supp}(\ ^{\on{loc}}\bV_{\bmu,(\hbar,\underline{c})}^{\Gamma_n})$ must be contained in exactly one component
$\fL_{\blambda,(\underline{c}^p-\hbar^{p-1}\underline{c})}^{(1)}$. It is easy to
 see that dimension of the support equals $n$ which is also the dimension of
$\fL_{\blambda,(\underline{c}^p-\hbar^{p-1}\underline{c})}^{(1)}$, thus
$$\on{supp}(\ ^{\on{loc}}\bV_{\bmu,(\hbar,\underline{c})}^{\Gamma_n})=
\fL_{\blambda,(\underline{c}^p-\hbar^{p-1}\underline{c})}^{(1)}.$$
It remains to prove that
$\blambda=\bmu$. It follows from the last displayed equality that
any finite dimensional  ${\mathbb G}_m^h$-equivariant quotient of
$\bV_{\bmu,(\hbar,\underline{c})}$ is supported at
the point $\blambda$ of the spectrum of the $p$-center of the Cherednik algebra
$\bH_{(\hbar,\underline{c})}$. On the other hand, $\bV_{\bmu,(0,\underline{c})}$
is the fiber over a prime over $p$ in $\on{Spec}\BZ[\sqrt[r]{1}]$
of the family of Verma modules over the family of Cherednik algebras.
Recall the family $\bM(\bmu)$ of baby Verma modules over the family of Cherednik
algebras over $\on{Spec}\BC\to\on{Spec}\BZ[\sqrt[r]{1}]$ introduced
in~\cite[4.1]{G0}.
This family of Cherednik algebras comes from one over
$\on{Spec}\BZ[\sqrt[r]{1}]$, and its fiber over a prime over $p$ coincides
with the family $\bH_{(0,\ul{c}),\sk}$ of Section~\ref{loca}.
The corresponding family of spherical Verma (resp. baby Verma) modules is
denoted $\bV_{\bmu,(0,\underline{c}),\sk}^{\Gamma_n}$ (resp.
$\bM(\bmu)_{\underline{c},\sk}^{\Gamma_n}$).
Now $\bV_{\bmu,(0,\underline{c}),\sk}^{\Gamma_n}$
surjects onto the spherical baby Verma module $\bM(\bmu)_{\underline{c},\sk}^{\Gamma_n}$.
The support of $\bM(\bmu)_{\underline{c},\sk}^{\Gamma_n}$ as a module over
the $\hbar=0$ spherical algebra $e\bH_{(0,\ul{c})}e=\sk[\bZ_{\ul{c}}]$ is the fixed point $\bmu$. Recall that the $p$-center of the spherical algebra
$e\bH_{(\hbar,\ul{c}),\sk}e$ equals $\sk[\bZ_{\underline{c}^p-\hbar^{p-1}\underline{c})}^{(1)}]$, and the
map from the spectrum of the $p$-center to the spectrum of the $\hbar=0$
spherical algebra is the Frobenius map~(Lemma~\ref{qHam}(c)), hence the $p$-support
of $\bM(\bmu)_{\underline{c},\sk}^{\Gamma_n}$ is also the fixed point $\bmu$.
 It follows
that $\blambda=\bmu$. This completes the proof of the lemma along with the
proof of the proposition. \qed

\proof of Proposition \ref{benz}. Recall the family of quiver varieties
$\varrho:\ \bY\to\BA^1_\hbar\times\BA^r_c$ introduced
in~Proposition~\ref{exist_deform}.
We need to show that
$$\fL_{\bmu}^{(1)}\subset \on{supp}(\ ^{\on{loc}}\sV_{\bmu}^{\Gamma_n})\subset\fL_{\bmu}^{(1)}\cup
\bigcup_{\blambda\precc\bmu}\fL_\blambda^{(1)}.$$

In view of Lemmas \ref{loseu} and \ref{benzrou} this would follow once we check that

i) the support of
$\on{supp}(\ ^{\on{loc}}\bV_{\bmu}^{\Gamma_n})$ is irreducible and

ii) $\on{supp}(\ ^{\on{loc}}\sV_{\bmu}^{\Gamma_n})= \on{supp}(\
^{\on{loc}}\bV_{\bmu}^{\Gamma_n})\cap\varrho^{-1}(1,0,\ldots,0)$.

Now, $^{\on{loc}}\bV_{\bmu}^{\Gamma_n}$ is an object in the derived
category of sheaves of modules over the Azumaya algebra $\bA$ concentrated in non-positive
homological degrees. Clearly, $^{\on{loc}}\sV_{\bmu}^{\Gamma_n}$
is obtained from  $^{\on{loc}}\bV_{\bmu}^{\Gamma_n}$ by derived
tensor product with skyscraper at $\hbar=1,\ \ul{c}=0$ over the ring of
 functions $\sk[\BA^1_\hbar\times\BA^r_c]$.
Lemma~\ref{exact} shows that the latter tensor product is concentrated in
 homological degree zero. It follows that in some  Zariski neighborhood of the
 special fiber $\varrho^{-1}(1,0,\ldots,0)$ the complex
$^{\on{loc}}\bV_{\bmu}^{\Gamma_n}$ is concentrated in homological degree zero and is torsion free as a module over the ring of functions
$\sk[\BA^1_\hbar\times\BA^r_c]$.
Thus every component of $\on{supp}(\ ^{\on{loc}}\bV_{\bmu}^{\Gamma_n})$
whose closure intersects $\varrho^{-1}(1,0,\ldots,0)$ maps
dominantly to $\BA^1_\hbar\times\BA^r_c$. Since the support is $\Gm$-invariant
under the $\Gm$-action of~Proposition~\ref{exist_deform}(e),
property (i) is established. Since $^{\on{loc}}\sV_{\bmu}^{\Gamma_n}= \
^{\on{loc}}\bV_{\bmu}^{\Gamma_n}\otimes_{\sk[\BA^1_\hbar\times\BA^r_c]}
\sk[\varrho^{-1}(1,0,\ldots,0)]$,
property (ii) is also clear.
 \qed

\section{The proof of the main Theorem}
\subsection{Procesi bundle}
\label{Proc}
Under the equivalence $\Phi_\sk$ of~(\ref{mult}) the structure sheaf
$\CO_{U^{(1)}}\rtimes\sk[\Gamma_n]\in D^b\Coh^{\Gm\times\Gamma_n}(U^{(1)})$ goes to
$\CE\in D^b\Coh^\Gm(Y^{(1)})$, to be called the {\em Procesi vector bundle}.
It splits into the direct sum $\CE=\bigoplus\CE_\blambda\otimes\rho_\blambda^*$
over the set of
$r$-partitions of $n$ numbering the irreducible $\sk[\Gamma_n]$-representations;
under the above equivalence $\CE_\blambda$ corresponds to
$\CO_{U^{(1)}}\otimes\rho_\blambda\in D^b\Coh^{\Gm\times\Gamma_n}(U^{(1)})$.

Recall that for an irreducible $\sk[\Gamma_n]$-module $\rho_\blambda$ we
denote by $\CF_\blambda$ the coherent sheaf
$\Phi_\sk(\CO_{L^{(1)}}\otimes\rho_\blambda)$.

We denote by $'L\subset U$ the $\Gamma_n$-invariant Lagrangian complement to
$L$ in $U$. Repeating the definitions of~Section~\ref{components} with $L$
replaced by $'L$ we obtain the Lagrangian components $'\fL_\bmu\subset Y$
numbered by the $r$-partitions of $n$. Note that the adjacency order on these
components is the inverse of the order in {\em loc. cit.}:
$'\fL_\blambda\cap\ '\fL_\bmu^\circ\ne\emptyset$ iff $\blambda\preceqg\bmu$.

We define the opposite Verma modules over $\sW_\sk\rtimes\sk[\Gamma_n]$ as
$'\sV_\bmu:=\sk[\ 'L]\otimes\rho_{\bmu}$. Under the Morita equivalence they
correspond to the $\sW_\sk^{\Gamma_n}$-modules
$'\sV_\bmu^{\Gamma_n}:=\sk[\ 'L]^{\rho_{\bmu}}$.
We have $\Ext^\bullet(\sV_\blambda,\ '\sV_\bmu)=0$ iff $\bmu\ne\ ^t\!\blambda$ 
(see \ref{components} for notations),
and $\Ext^n(\sV_\blambda,\ '\sV_{^t\!\blambda})\simeq\sk$, since
$\Ext_{\sk[U]}^\bullet(\sk[L],\sk[\ 'L])$ is a one-dimensional $\sk$-vector
space in degree $n$ where $\Gamma_n$ acts  by  the  representation
$\det L'$, and $\rho_\bmu\otimes\det L\simeq\rho_{^t\!\bmu}$.

Also we  define the sheaf
$'\CF_\bmu:=\Phi_\sk(\CO_{'L^{(1)}}\otimes\rho_{\!\bmu})$.

\begin{prop}
\label{renamed}
We have $\on{supp}\CF_\blambda=\on{supp}(\ ^{\on{loc}}\sV_\blambda)$.
\end{prop}

\proofpt Recall that the $\Gamma_n$ equivariant splitting bundle
$\CE_{orb}$ used in the definition of the equivalence $\Phi_\sk$ is assumed
to satisfy a compatibility stated after the proof of Lemma \ref{lem_comm} (end
of section \ref{sec22}).
It is clear that a Verma module $\sV_\blambda$ is isomorphic to
$\CO(L)\otimes _{\CO_{L^{(1)} } } (\CO_{L^{(1)}}\otimes\rho_{\blambda})$.
Thus the Proposition follows from the commutative diagram in Lemma \ref{lem_comm}(a).
%
\qed

\begin{cor}
\label{benzroukao}
If the Lagrangian component $\fL^{(1)}_\bmu$ lies in the support of
$\CF_\blambda$, then $\bmu\preceqc\blambda$. If the Lagrangian component
$'\fL^{(1)}_\bmu$ lies in the support of
$'\CF_\blambda$, then $\bmu\succeqc\ ^t\blambda$.
\end{cor}

\proofpt
We  argued in the proof
of~Proposition~\ref{benz} that the support of the localized Verma module
$^{\on{loc}}\sV_\blambda$ consists of $\fL^{(1)}_\blambda$ and some smaller
Lagrangian components $\fL^{(1)}_\bmu,\ \bmu\precc\blambda$. Thus Corollary
follows from Proposition \ref{renamed}.

The second part is proved similarly.
\qed

\subsection{Proof of Theorem~\ref{unique}}
\label{210}
Note that the composed projection
$Y\to\BA^{2n}/\Gamma_n\to\BA^n/\Gamma_n$ is flat since all the
fibers are of the same dimension $n$. It follows that
$\overline\CE_{\rho_\bmu}\cong \CF_\bmu$ and $\overline\CE{}'_{\rho_\bmu}\cong\ '\CF_\bmu$;
the required properties of $\CF_\bmu$, $\ '\CF_\bmu$ have been established
in~Corollary~\ref{benzroukao}. \qed

\subsection{Wreath Macdonald polynomials}
\label{wreath}
The equivalence $\Phi_\sk$ extends to the same named equivalence
$D^b\Coh^T(U^{(1)})\iso D^b\Coh^T(Y^{(1)})$ where $T$ is the 2-dimensional
torus with coordinates $(q,t)$ acting on $\BA^2$ such that $\Gm$ is the
diagonal subtorus $\Gm=\{(q,q)\}$, while
${\mathbb G}^h_m$ is the antidiagonal subtorus $\{(q,q^{-1})\}$.
We keep the same name for the
induced isomorphism of the $K$-groups of the categories in question.
We identify the equivariant $K$-group $K^T(pt)$ with $\BZ[q^{\pm1},t^{\pm1}]$.
We denote by $\iota_\blambda^*:\ K^T(Y^{(1)})\to K^T(pt)$ the fiber at the
torus fixed point $\blambda\in Y^{(1)}$.

Note that $\iota_\bmu^*\CE$ carries an action of $\Gamma_n$ since
$\CE$ carries a fiberwise action of $\Gamma_n$ by construction. Thus
$[\iota_\bmu^*\CE]\in\BZ[q^{\pm1},t^{\pm1}][K(\Gamma_n)]$.

\begin{cor}
\label{rom}
$\iota_\bmu^*[\CE]\otimes\sum_i(-q)^i[\Lambda^i(\ 'L^{(1)})]\in
\BZ[q^{\pm1},t^{\pm1}][\rho_\blambda,\ ^t\blambda\succeqc\bmu]$; also,\\
$\iota_\bmu^*[\CE]\otimes\sum_i(-t)^i[\Lambda^i(L^{(1)})]\in
\BZ[q^{\pm1},t^{\pm1}][\rho_\blambda,\ ^t\blambda\succeqc\ ^t\bmu]$.
\end{cor}

\proof (cf.~\cite[Proof of~Proposition~5.4.2]{H}).
The Koszul resolution shows that
$[\CO_{'L^{(1)}}]=[\CO_{U^{(1)}}]\otimes\sum_i(-t)^i[\Lambda^i(L^{(1)})]$.
We have $\operatorname{Tor}^{\CO_{L^{(1)}}}(\iota^*_\bmu(\CE),\sk)=
\operatorname{Tor}^{\CO_{U^{(1)}}}(\iota^*_\bmu(\CE\otimes\CO_{'L^{(1)}}),\sk)$.
By the local duality, the second statement of the corollary is equivalent to
$\operatorname{Ext}_{\CO_{U^{(1)}}\rtimes\sk[\Gamma_n]}(\CO_{'L^{(1)}}\otimes
\rho_\blambda,\iota^*_\bmu(\CE))=0$ unless $\blambda\succeqc\ ^t\bmu\
\Leftrightarrow\ ^t\blambda\preceqc\bmu$. Since $\Phi_\sk(\iota^*_\bmu(\CE))$ is
the skyscraper sheaf at the point $\bmu\in Y$, the latter Ext-vanishing
follows from $\operatorname{supp}\Phi_\sk(\CO_{'L^{(1)}}\otimes\rho_\blambda)
\subset\bigcup_{\bmu\succeqc\ ^t\blambda}\ '\fL_\bmu$.

The first statement of the corollary is proved similarly.
\qed

\section{Errata to \cite{BK}}
\label{erra}
We use this opportunity to correct some mistakes in \cite{BK}.

Proof of Lemma 2.12 on p. 8: the action of ${\mathbb G}_m$ on $X_R$ induces
an action of the formal completion of ${\mathbb G}_m$ on the formal completion
${\mathcal X}$, cf. Appendix
(section \ref{volgd}). So references to ${\mathbb G}_m$-equivariant structures
on the bundle on $\mathcal X$ in the proof should be replaced by ones
to equivariance with respect to the formal completion of ${\mathbb G}_m$.

Proposition 3.8 on p. 12: ``restriction to the generic fiber of the formal scheme $\hat{X}$"
should be removed, the statement of the Proposition should read as follows:
"The localization $\O_\hbar\otimes_{\sk[[\hbar]]}\sk((\hbar)) $ is a sheaf of Azumaya algebras over
the sheaf of commutative rings $\O_{\hat{X}}\otimes _{\sk[[\hbar]]}\sk((\hbar))$."

The most essential correction\footnote{We thank Ivan Losev for pointing this out.}
 is for section 5.1.3 on p. 20 which provides the central step
in the proof of Proposition 5.3. It is claimed  there that the  quantization
of Kleinian singularity whose existence is claimed in  Proposition 5.3
is obtained from the Weyl algebra of the vector space $M$ by Hamiltionian reduction,
using the quantum moment map described in Example 5.7. This is false; e.g. in the example
$\Gamma=\{\pm 1\}$ acting on the two dimensional symplectic vector space, (the fiber at $\hbar=1$ of)  the quantization
arising from Example 5.7 is identified with  the enveloping algebra
$U(\mathfrak{sl}(2,\sk))$ reduced
at the singular Harish-Chandra central character $-\rho=-1$, while
$\Gamma$-invariants
in the Weyl algebra are isomorphic to the reduction of
$U(\mathfrak{sl}(2,\sk))$ by the regular central
character $-\frac{1}{2}$.
The argument breaks down as the displayed isomorphism
${\mathcal R}_\hbar \cong \O_\hbar\otimes_{\O_{X_0}^p}\otimes \O_{V_0}^p$ does not hold,
because, contrary to the mistaken statement in \cite{BK},
the commutative ring $\O(V^{(1)})[[\hbar]]$ does not act on $\CR_\hbar$ in  a natural way.

In order to get the quantization claimed in Proposition 5.3 one has to use a different
quantum moment map, which differs from the one in \cite[Example 5.7]{BK}
by adding the character
 described in the next paragraph, this fact goes back to~\cite{CBH}\footnote{Also the following special case of this fact
coincides with the special case $n=1$ of Lemma \ref{qHam}: we consider  the  quiver of type $\widehat{A_{r-1}}$, so $I=\BZ/r\BZ$, $i_0=0$, the orientation is chosen so that $i$ is connected to
$i+1$ by an oriented edge for $i\in I$; $\delta_i=1\ \forall i\in I$, thus
$\partial_i=0\ \forall i\in I$.} in the case of characteristic zero, the case of large positive characteristic
follows by the argument in the proof of Lemma \ref{qHam}.

Consider an extended Dynkin quiver with the set of vertices $I$ and the
extending vertex $i_0\in I$.  We set $\delta_i,\ i\in I$, to be the $i$-th coordinate of the
corresponding minimal imaginary root vector (e.g. $\delta_{i_0}=1$). We choose an orientation
of the extended Dynkin quiver. We define the {\em defect vector}
$\partial$ with coordinates $\partial_i$ as follows:
$\partial_i:=-\delta_i+\sum\delta_j$ where the sum is over the set of vertices
$j$ such that there is an arrow $i\to j$.  Finally, $\Gamma\subset SL(2)$ is the finite
subgroup corresponding to our quiver. Then the character producing the
hamiltonian reduction equal to the $\Gamma$-invariants in the first Weyl
algebra has the following components: then
$\chi_{i_0}=\delta_{i_0}/|\Gamma|-\partial_0-1,\
\chi_i=\delta_i/|\Gamma|-\partial_i\ (i\ne i_0)$.
 
\section{Appendix: Equivariant structure on rigid vector bundles}
\label{volgd}

\centerline{\bf {by Vadim Vologodsky}}
\vskip 6pt

Let $\TTT=(\Gm)^d$ be a $d$-dimensional torus.
We begin with the following result.

\begin{lem}
\label{odin}
Let $X$ be a proper (not necessarily smooth) scheme over an algebraicaly
closed field $\sk$, $E$ a vector bundle over $X$ such
that $\Ext^1(E;E)=0$. Assume that $X$ is endowed with a
$\TTT$-action  $m: \TTT \times X \to X$.
Then the action of  $\TTT$ lifts to $E$, i.e.
$E$ has a $\TTT$-equivariant structure.
\end{lem}
\proofpt First, we claim that $E$ is $\TTT$-invariant, that is
$m^* E \simeq p_2^* E$  locally with respect to $\TTT$.

Indeed, since $\Ext^1(E;E)=0$, the claim is true if we replace
$\TTT$ by its formal completion at $1$.
Now, let us consider the vector bundle $F= \Hom(m^* E, p_2^* E)$.
Since $X$ is proper, $p_{1 *} F$ is a coherent sheaf on  $\TTT$.
Theorem on Formal Functions implies that the
formal completion of $p_{1 *}F$
 at $1$ is a vector bundle. Therefore $p_{1 *}F$ is a vector bundle over
some
open neighborhood $U$ of $1$ in  $\TTT$. Further, shrinking
$U$, we may assume that there exists a section $s$ of   $p_{1 *}F$ over
$U$ which is equal to $\Id \in (p_{1 *}F)_1$ at the point $1\in  \TTT$.
It implies that   $m^* E \simeq p_2^* E$ over some open  neighborhood of
$1$. Since the subgroup of $\TTT$ generated by an open neighborhood
of $1$ coincides with $\TTT$, the claim follows.

We conclude that there exists an extension of affine algebraic groups
$$1\to \Aut(E) \to G \to  \TTT \to 1$$
such that the action of $G$ on $X$ (through  $\TTT$) lifts to $E$.

Note that the group $\Aut(E)$ is smooth and connected (since it is an
open subscheme
of the  affine space $\End(E)$). Hence the lemma follows from the next one.

\begin{lem}\label{dva}
For a smooth  connected group $H$,
any extension
  $$1\to H \to G \stackrel{h}{\to}  \TTT \to 1$$
splits in the weak sense, i.e. there exists
a homomorphism $i:\TTT \to   G$ such that $h\circ i = \Id$.
\end{lem}
\proofpt The result is well known. We include the argument for the
readers' convenience.

Any affine algebraic group is a semidirect product of a reductive group
and a unipotent one. Therefore it is enough to prove the lemma in the following
two cases.
\begin{enumerate}
\item{ $H$ is a unipotent  group.
In this case the lemma follows from the structural result cited above:
if $G= U\times \TTT$ is the semidirect factorization, the morphism
$h$ factors through the second factor and, since the kernel of $h$ is
connected, there exists a splitting.}

\item{$H$ is  a reductive group.
In this case, $\rk(G)=\rk(H)+\dim (\TTT)$. Therefore, the intersection of a
maximal torus  $T\subset G$ with $H$ is a toral subgroup of $H$
of dimension $\rk(H)$. Hence its  connected component is a Cartan
subgroup of $H$;
the intersection $T\cap H$ is contained in the centralizer of that
Cartan subgroup,
hence it coincides with the Cartan subgroup. Thus, the lemma reduces
 to the obvious case when $G$ is a torus. }
 \end{enumerate}
\qed

Let
${\mathcal X}\to{\mathcal Y}$ be a morphism of schemes of finite type
over an algebraically closed field.  Let $y_l$ (resp. $\hat y$) be
the $l$-th infinitesimal neighborhood (resp. the formal neighborhood)
of a closed  point $y\in{\mathcal Y}$, and let  ${\mathcal X}_{y_l}$
(resp. ${\mathcal X}_{\hat y}$)
be the preimage of $y_l$ (resp. $\hat y$)  in $\mathcal X$.  Assume
that  ${\mathcal X}_{y_l}$, ($l=1,2,\cdots$), are endowed with
compatible $\TTT$-actions
$m_l:\ \TTT\times{\mathcal X}_{y_l}\to{\mathcal X}_{y_l}$\footnote{In
particular,  ${\mathcal X}_{\hat y}$ is equipped with an action of the
formal completion $\hat  \TTT$ of $  \TTT$.}.
We will need a version of the notion of an equivariant bundle for
formal neighborhoods. By a $\TTT$-equivariant  structure  on
a vector bundle $E$ over  ${\mathcal X}_{\hat y}$
we mean a compatible system of $\TTT$-equivariant structures on $E_l:=
E_{|{\mathcal X}_{y_l}}$. Equivalently,  this is  an action of
$\hat  \TTT$ on $E$ lifting the $\hat  \TTT$-action on
 ${\mathcal X}_{\hat y}$ such that for each $l$ the $\hat
\TTT$-action on $E_l$ can be extended to an action of $ \TTT$.
\begin{prop}
\label{vologod}
Assume that
${\mathcal X}\to{\mathcal Y}$ is proper and that $\Ext^1(E,E)=0$. Then
$E$ admits a  $\TTT$-equivariant structure.
\end{prop}

\proofpt
Consider the projective system of algebraic groups
$H_l=\Aut(E_l)$. Since
$\End_\sk(E_l)$ is a finite-dimensional $\sk$-vector space,
the projective system above
satisfies the Mittag-Leffler condition, {\it i.e.},
for any integer $l$ there exists an integer $N\geq l$ such that
$$\im( H_k \to H_l) =\im( H_N \to H_l), \quad\forall k>N.$$
Similarly, using that $\Ext^1(E,E)=0$ and  that $\Ext^1(E_l, E_l)$ is
finite-dimensional, we conclude that
$$\im(\Ext^1(E_N, E_N)\to\Ext^1(E_l, E_l))=0,$$
for $N$ large enough.

Using the above assertions and
arguing as in the proof of Lemma ~\ref{odin},
we conclude that, for every $l$, the bundle $E_l$ admits a
$\TTT$-equivariant structure.

Next, consider a projective system of extensions
$$1\to H_l \to G_l \to \TTT\to1,$$
where,  for a scheme $S$, $G_l(S)$ is the group of pairs $(t: S \to
\TTT, \theta)$, where $\theta$ is an isomorphism between $(m_l \circ
(t\times Id))^* E_l $ and the pullback of $E_l$ on $S\times  {\mathcal
X}_{y_l}$ with respect to the  second projection.
 The Mittag-Leffler property of $H_l$ together with Lemma ~\ref{dva}
ensures the existence of a weak splitting $\TTT\to \varprojlim G_l$
which defines a $\TTT$-equivariant structure on $E$.
\qed

\bigskip

\footnotesize{ {\bf R.B.}: Department of Mathematics, Massachusetts Institute
of Technology,\\ Cambridge MA 02139, USA and\\
National Research University Higher School of Economics,
International Laboratory of Representation
Theory and Mathematical Physics,
20 Myasnitskaya Ulitsa, Moscow 101000, Russia;\\ 
{\tt bezrukav@math.mit.edu}}

\footnotesize{ {\bf M.F.}: IMU, IITP, and
National Research University Higher School of Economics,\\ Department of
Mathematics, 20 Myasnitskaya st, Moscow 101000 Russia;\\ 
{\tt fnklberg@gmail.com}}


\footnotesize{ {\bf V.V.}: Department of Mathematics, Oregon University,
Eugene OR 97403, USA;\\ 
{\tt vvologod@uoregon.edu}}

\end{document}